\documentclass{amsart}
\usepackage{mathrsfs}
\usepackage{amsfonts}
\usepackage{graphicx}
\usepackage{amsmath}
\usepackage{amssymb}
\usepackage{bm}
\usepackage{color}
\usepackage{mathdots}
\usepackage[scr=boondoxo]{mathalfa}

\usepackage[linktocpage=true
]{hyperref}
\hypersetup{colorlinks,linkcolor=blue,urlcolor=cyan,citecolor=blue,
pdftitle={Twisted Yangians of quasi-split type DI and Drinfeld type current relations},
pdfauthor={Vidas Regelskis}}

\newtheorem{theorem}{Theorem}[section]
\newtheorem{thm}[theorem]{Theorem}

\newtheorem{crl}[theorem]{Corollary}

\newtheorem{lemma}[theorem]{Lemma}
\newtheorem{prop}[theorem]{Proposition}

\numberwithin{equation}{section}
	
\theoremstyle{remark}
\newtheorem{remark}[theorem]{Remark}

\newcommand{\nc}{\newcommand}

\nc{\ot}{\otimes}
\nc{\op}{\oplus}
\nc{\ol}{\overline}

\nc{\mc}{\mathcal}
\nc{\ms}{\mathsf}
\nc{\mf}{\mathfrak}
\nc{\mb}{\mathbf}
\nc{\bb}{\mathbb}
\nc{\mr}{\mathscr}

\nc{\al}{\alpha}
\nc{\bet}{\beta}
\nc{\eps}{\epsilon}
\nc{\del}{\delta}
\nc{\ga}{\gamma}
\nc{\Ga}{\Gamma}
\nc{\ka}{\kappa}
\nc{\vka}{\varkappa}
\nc{\la}{\lambda}
\nc{\om}{\omega}
\nc{\si}{\sigma}
\nc{\Ups}{\upsilon}

\nc{\ud}{\underline}
\nc{\tl}{\tilde}
\nc{\wt}{\widetilde}
\nc{\wh}{\widehat}
\nc{\pd}{\partial}

\nc{\C}{\mathbb{C}}
\nc{\N}{\mathbb{N}}
\nc{\Z}{\mathbb{Z}}

\nc\mfgl{\mathfrak{gl}}
\nc\mfsl{\mathfrak{sl}}
\nc\mfsp{\mathfrak{sp}}
\nc\mfo{\mathfrak{o}}

\nc{\lan}{\langle}
\nc{\ran}{\rangle}
\nc{\tw}{{\sf tw}}

\nc{\equ}[1]{\begin{equation}#1\end{equation}}
\nc{\eqa}[1]{\begin{equation}\begin{alignedat}{50}#1\end{alignedat}\end{equation}}
\nc{\eqn}[1]{\begin{equation*}\begin{alignedat}{50}#1\end{alignedat}\end{equation*}}
\nc{\eqg}[1]{\begin{equation}\begin{gathered}#1\end{gathered}\end{equation}}

\nc{\ali}[1]{\begin{alignat}{50}#1\end{alignat}}
\nc{\als}[1]{\begin{subequations}\begin{alignat}{50}#1\end{alignat}\end{subequations}}
\nc{\aln}[1]{\begin{alignat*}{50}#1\end{alignat*}}

\nc{\gat}[1]{\begin{gather}#1\end{gather}}
\nc{\gas}[1]{\begin{subequations}\begin{gather}#1\end{gather}\end{subequations}}
\nc{\gan}[1]{\begin{gather*}#1\end{gather*}}

\nc\el{\nonumber\\}
\nc\nn{\nonumber}

\nc{\qu}{\quad}
\nc{\qq}{\qquad}

\nc{\iso}{\stackrel{\sim}{\longrightarrow}}
\nc{\into}{\hookrightarrow}
\nc{\onto}{\twoheadrightarrow}

\DeclareMathOperator{\tr}{tr}
\DeclareMathOperator{\Sym}{Sym}
\DeclareMathOperator{\id}{id}
\DeclareMathOperator{\gr}{gr}
\DeclareMathOperator{\diag}{diag}
\DeclareMathOperator{\End}{End}

\renewcommand{\,}{\kern 0.1em} %
\begin{document}

\title[Twisted Yangians of quasi-split type DI]{Twisted Yangians of quasi-split type DI \\ and Drinfeld type current relations}

\author{Vidas Regelskis}
\address{Department of Physics, Astronomy and Mathematics, University of Hertfordshire,
Hatfield AL10 9AB, UK and
Institute of Theoretical Physics and Astronomy, Vilnius University, Saul\.etekio av.~3, Vilnius 10257, Lithuania.}
\email{vidas.regelskis@gmail.com}


\begin{abstract}
We study the twisted Yangian associated with the quasi-split symmetric pair $(\mfo_{2n}, \mfo_{n-1}\op\mfo_{n+1})$ of type DI. Working in the transposed $R$-matrix presentation of \cite{Re26}, we show that the generating matrix $S(u)$ admits no Gaussian decomposition. We decompose instead the matrix $\wh{S}(u) = G\,S(u)$, where $G$ is the constant term of $S(u)$. A chain of quasi-determinantal embeddings lowers the rank within this family of pairs and ends at $(\mfo_4,\mfo_3)$. Together with a direct study of a degenerate block of type AI, this gives a presentation of the special twisted Yangian $SY^\tw(\mfo_{2n},\mfo_{n-1}\op\mfo_{n+1})$ by Drinfeld type current relations. The two fork nodes of the Dynkin diagram carry Drinfeld's presentation of $Y(\mfsl_2)$. All other relations are those of split type DI, except that the Serre relation with a doubled fork node has zero right hand side. For $n\ge3$ we obtain an explicit isomorphism between the $R$-matrix presentation of \cite{GR16} and the Drinfeld presentation of Lu and Zhang \cite{LZ25}. The Drinfeld generators are expressed through the Gaussian series of $\wh{S}(u)$. With quasi-split type AIII \cite{LZ25}, this completes the explicit comparison for the classical quasi-split types with a nontrivial diagram involution.
\end{abstract}

\maketitle
\setcounter{tocdepth}{1}
\tableofcontents
\setcounter{tocdepth}{2}

\setlength{\medmuskip}{2.0mu plus 1.0mu minus 1.0mu}


\section{Introduction} \label{sec:intro}

Twisted Yangians are coideal subalgebras of Yangians attached to symmetric pairs. Their
$R$-matrix presentation builds on Sklyanin's reflection equation \cite{Sk88}. In this
presentation they were introduced by Olshanski \cite{Ol92} for the symmetric pairs of types AI
and AII (see also \cite{MNO96}), by Molev and Ragoucy \cite{MR02} for type AIII, and by Guay and
the author \cite{GR16} for the remaining classical types B, C, D. The cases of $\mfsp_2$, $\mfo_3$ and
$\mfo_4$, except the diagonal pair $(\mfo_4,\mfo_3)$, were studied in \cite{GRW16}.

Drinfeld type current presentations of twisted Yangians appeared only recently. They are
modelled on Drinfeld's new realisation of Yangians \cite{Dr87}, with Cartan and root currents
attached to the nodes of the Satake diagram. Lu, Wang and
Zhang obtained one for type AI \cite{LWZ25a} and, by degenerating the affine $\imath$quantum
groups of \cite{LW21, Zh22}, formulated one for all split types \cite{LWZ25b}. Lu and Zhang then
formulated one for the quasi-split types, and proved its equivalence with the $R$-matrix presentation
only in type AIII \cite{LZ25}. Minimalistic presentations and the coideal structure were studied
by Lu \cite{Lu26a} in the split types and by Huang and Lu \cite{HL26} in the quasi-split types.

In the companion paper \cite{Re26} we introduced the \emph{transposed presentation} of the
twisted Yangians of types B, C, D. It is governed by a twisted reflection equation
rather than by the (non-twisted) reflection equation of \cite{GR16}. There we derived the current relations for the split pairs of
types BI, CI, DI by Gaussian decomposition of the generating matrix $S(u)$, as in
\cite{JLM18, LWZ25a}. The constant term $G$ is diagonal and
invertible in the split case. Hence every leading principal submatrix of $S(u)$ is invertible.

This paper studies the twisted Yangian of the symmetric pair
\equ{
(\mfo_{2n},\, \mfo_{n-1}\op\mfo_{n+1}) , \qq n \ge 2 ,
\label{our-pair}
}
of type DI. It is quasi-split but not split: its Satake diagram is the Dynkin diagram of $D_n$
with all nodes white, and the diagram involution $\tau$ interchanges the two fork nodes. For
$n\ge4$ it is of type DII$_n$ in \cite{LZ25}, and for $n=3$ of type AIII$_3$. For $n=2$ it is the diagonal pair
$(\mfo_4,\mfo_3)\cong(\mfsl_2\op\mfsl_2,\Delta\mfsl_2)$. Twisted Yangians of diagonal pairs are
the achiral twisted Yangians of MacKay and the author \cite{MaRe11}.

For the pair \eqref{our-pair} the matrix $G$ is the permutation matrix of $(n,n+1)$, and $S(u)$
has no Gaussian decomposition in any ordering of the basis. We decompose instead $\wh{S}(u) := G\,S(u)$,
whose constant term is $G^2=I$. This matrix is compatible with the quasi-determinantal
embeddings of \cite{Re26}. The embeddings stay inside the family \eqref{our-pair} and end at the
diagonal pair, for which $SY^\tw(\mfo_4,\mfo_3)\cong Y(\mfsl_2)$ (Theorem \ref{T:o4o3}). For
$n\ge4$, the embeddings reduce the relations between the fork nodes and their neighbour to $n=3$. At
$n=3$, we use the leading $3\times3$ block of $S(u)$. It satisfies the twisted reflection equation of
type AI. The split DI case also uses a block of type AI \cite[\S5.1 and \S5.3]{Re26}. Here the block
is degenerate: its last diagonal Gaussian series vanishes. So the relations of the split case cannot be imported, and we derive them directly from
the twisted reflection equation (Section~\ref{sec:deg}).

The main result, Theorem \ref{T:iso}, presents
$SY^\tw(\mfo_{2n},\mfo_{n-1}\op\mfo_{n+1})$ by Drinfeld type current relations. The two fork nodes carry Drinfeld's presentation of $Y(\mfsl_2)$. All other pairs of nodes
satisfy the relations of split type DI \cite[Thm.~5.14]{Re26}, with one exception: the Serre
relation with a doubled fork node has zero right hand side. For $n\ge3$, Corollary \ref{C:LWZ}
identifies the twisted Yangian with the algebra $\mc{Y}^\imath$ of Lu and Zhang \cite{LZ25}, by an
explicit isomorphism on Drinfeld generators. With \cite[Thm.~B]{LZ25} for type AIII, this relates
the $R$-matrix and Drinfeld presentations explicitly in every classical quasi-split type with
$\tau\ne\id$. For $n\ge3$, Corollary \ref{C:iso-ext} gives the corresponding presentation of the extended twisted Yangian.

While this paper was being completed, Huang and Lu \cite{HL26} embedded $\mc{Y}^\imath$ into the
Yangian and identified it with a twisted Yangian in the $J$ presentation of MacKay \cite{Mac02}.
They start from the Drinfeld presentation. We start from the $R$-matrix presentation and derive
the Drinfeld one. Together the two results give an isomorphism between the $R$-matrix and $J$
presentations for $n\ge3$ (Corollary \ref{C:J}). The classification of twisted Yangians by Lu \cite{Lu26b} gives an
alternative proof of this (Remark \ref{R:Lu}).

This paper is organised as follows. Section~\ref{sec:setup} sets up the pair \eqref{our-pair},
the matrix $\wh{S}(u)$, the centre and the reduction to lower rank. Section~\ref{sec:low} treats
the low rank case $(\mfo_4,\mfo_3)$, and Section~\ref{sec:deg} the degenerate block.
Section~\ref{sec:Dr} proves the presentation and compares it with \cite{LZ25} and \cite{HL26}. It
also treats the extended twisted Yangian, the fundamental representation and the coideal
structure.


\section{The quasi-split pair and the matrix \texorpdfstring{$\wh{S}(u)$}{S(u)}} \label{sec:setup}


\subsection{Notation} \label{sec:notation}

We follow the conventions of \cite[\S2]{Re26}, to which we refer for all statements recalled in
this subsection. Throughout, $N = 2n$ with $n\ge2$, the Lie algebra is $\mfo_{2n}$, so that
$\theta = 1$, and the matrix $G$ of \eqref{G} below satisfies $G^2 = I$, so that the sign
$(\pm)$ of \cite[\S2.2]{Re26} equals $1$. Both are suppressed in the displays that follow. 
We write $E_{ij}$ for the matrix units of $\End(\C^N)$, $e_i$ for the standard basis vectors,
$\ol{\imath} := N-i+1$, and
\equ{
I := \sum_{i,j} E_{ii}\ot E_{jj} , \qu
P := \sum_{i,j} E_{ij}\ot E_{ji} , \qu
Q := \sum_{i,j} E_{\ol{\jmath}\,\ol{\imath}} \ot E_{ji} , \qu
J := \sum_{i} E_{i\hspace{0.03em}\ol{\imath}} ,
}
the sums being over $1\le i,j\le N$. Set
\equ{
R(u) := I - u^{-1}P - (\ka-u)^{-1}Q
}
with $\ka := N/2 - \theta = n-1$ and write $R'(u) = R^{t_1}(u) = R^{t_2}(u)$ for its partial
transpose, $t$ denoting the ordinary matrix transpose.

Given a matrix $G\in\End(\C^N)$ as in \cite[\S2.2]{Re26}, the extended twisted Yangian
$X^\tw(\mfo_{2n},G)$ is the unital associative $\C$-algebra with generators $s^{(r)}_{ij}$,
$1\le i,j\le N$ and $r\in\N$, assembled into
\equ{
S(u) := \sum_{1\le i,j\le N} E_{ij}\ot s_{ij}(u) , \qq
s_{ij}(u) = g_{ij} + \sum_{r\ge1} s^{(r)}_{ij}u^{-r} ,
\label{S(u)}
}
subject to the twisted reflection equation %
\equ{
R(u-v)\,S_1(u)\,R'(\ka-u-v)\,S_2(v) = S_2(v)\,R'(\ka-u-v)\,S_1(u)\,R(u-v) .
\label{RE}
}
We write $S(\infty)$ for the constant part of $S(u)$, so that $S(\infty) = G$ by \eqref{S(u)}.
The expanded form of \eqref{RE} in terms of the series $s_{ij}(u)$ is
\cite[(2.19)]{Re26}, which we shall refer to freely. The special twisted Yangian
$SY^\tw(\mfo_{2n},G)$ is the quotient of $X^\tw(\mfo_{2n},G)$ by the two-sided ideal generated by the
symmetry relation, in which $G(u)$ is the matrix \eqref{G(u)} below, %
\ali{ \label{sym}
S^t(u) = S(\ka-u) + \frac{S(u)-S(\ka-u)}{2u-\ka} + \frac{\tr(G(u)J)\,S(\ka-u) - \tr(S(u)J)J }{2u-2\ka} 
}
and the unitarity relation
\equ{
S(u)\,J\,S(-u)\,J = I .
\label{unit}
}
The mapping %
\equ{
S(u) \;\mapsto\; T(u-\ka/2)\,G(u)\,T^t(\ka/2-u)
\label{S->TGT}
}
embeds $SY^\tw(\mfo_{2n},G)$ into the Yangian $Y(\mfo_{2n})$, where $T(u)$ is the (normalised)
generating matrix of $Y(\mfo_{2n})$ and $G(u)$ is the matrix-valued solution of the twisted reflection equation
recalled in \eqref{G(u)} below. The image satisfies the unitarity relation \eqref{unit} since
$T(u)$ is normalised and $G(u)\,J\,G(-u)\,J=I$. If $T(u)$ is instead the (non-normalised) generating matrix of the
extended Yangian $X(\mfo_{2n})$, the same formula gives an isomorphism of
$X^\tw(\mfo_{2n},G)/(\ms{c}(u)-1)$ onto its image \cite[Thm.~4.2, Cor.~5.2]{GR16}, where $\ms{c}(u)$ is
the central series of Section~\ref{sec:centre}. We refer to
both maps by \eqref{S->TGT}; the matrix used shows which one is meant. Write $F_{ij} := E_{ij} - E_{\ol{\jmath}\,\ol{\imath}}$ for the standard generators of
$\mfo_{2n}$, satisfying $F_{ij} + F_{\ol{\jmath}\,\ol{\imath}} = 0$.

Both algebras carry the exhaustive multiplicative filtration determined by
$\deg s^{(r)}_{ij} := r-1$, with $\mr{F}_{-1}=0$. For $A\in\mr{F}_k$, write $\bar A$ for its
symbol in $\mr{F}_k/\mr{F}_{k-1}$.

We recall that the mapping
\equ{
\bar{s}^{(r)}_{ij} \mapsto F^{(\rho,r-1)}_{ij} + \ga^{(r-1)}_{ij} , \qu
F^{(\rho,r)}_{ij} := \sum_{1\le a\le N}\big(F_{ia}g_{aj} - (-1)^r g_{ia}F_{ja}\big)x^r , \qu
\ga^{(r)}_{ij} := \del_{r0}\tfrac{q-p}{4}\big(g_{ij}-\del_{i\ol{\jmath}}\big) ,
\label{SY-graded-iso}
}
defines an isomorphism $\gr SY^\tw(\mfo_{2n},G) \cong U(\mfo_{2n}[x]^\rho)$ of the associated graded
algebra with respect to this filtration. Here $\mfo_{2n}[x]^\rho$ is the Lie subalgebra of
$\mfo_{2n}[x]$ spanned by the elements $F^{(\rho,r)}_{ij}$, that is, the fixed points of the
involution $X x^r\mapsto(-1)^r\rho(X)\,x^r$, with $\rho$ and $(p,q)$ as in \eqref{pair} below. The term $\ga^{(r)}_{ij}$ is
visible only when $p\ne q$, and its normalisation is recorded in \eqref{S1} below. The mapping
\equ{
\bar{s}^{(r)}_{ij} \mapsto F^{(\rho,r-1)}_{ij} + \tfrac12\,g_{ij}\,\zeta_{r-1}  + \ga^{(r-1)}_{ij}
\label{X-graded-iso}
}
defines an isomorphism
$\gr X^\tw(\mfo_{2n},G) \cong U(\mfo_{2n}[x]^\rho)\ot\C[\zeta_0,\zeta_1,\ldots]$, where
$\C[\zeta_0,\zeta_1,\ldots]$ is a polynomial algebra in countably many indeterminates with
$\deg\zeta_r = r$. Both statements are \cite[(2.29), (2.30)]{Re26}, which rest on \cite{GR16}. %

Last, we recall the symmetrising bracket and the symmetrisation over spectral
parameters:
\[
\{f(u)\}^u := f(u)+f(-u) , \qq
\Sym_{v,w} f(v,w) := f(v,w) + f(w,v) .
\]


\subsection{The symmetric pair} \label{sec:pair}

We take $(p,q) = (n-1,n+1)$ in the type DI entry of \cite[\S2.2]{Re26}, so that the symmetric
pair is
\equ{
(\mfo_{2n},\, \mfo_{n-1}\op\mfo_{n+1}) , \qq
\rho(X) = -G\,X^t\,G^{-1} , \qq \mfo^\rho_{2n} = \mfo_{n-1}\op\mfo_{n+1} ,
\label{pair}
}
with
\equ{
G := \sum_{i=1}^{n-1}\big(E_{ii} + E_{\ol{\imath}\ol{\imath}}\big) \;+\; E_{n,n+1} + E_{n+1,n} .
\label{G}
}
Thus $G$ is the permutation matrix of the transposition $(n,n+1)$; note that $\ol{n} = n+1$. In
particular
\equ{
G^t = G , \qq G^2 = I , \qq GJ = JG , \qq \tr(GJ) = 2 , \qq \tr(J^2)=N .
\label{G-props}
}
For $n=2$ we have $\mfo_4 \cong \mfsl_2\op\mfsl_2$ and $\rho$ is conjugate to the flip of the
two summands, so that \eqref{pair} is the diagonal pair $(\mfsl_2\op\mfsl_2,\Delta\mfsl_2)$. This
is the only case deferred in \cite[\S4.1]{GRW16}; it is considered in Section~\ref{sec:low}.

The associated matrix-valued solution of the twisted reflection equation is
\equ{
G(u) := \frac{(p-q)J + 4u\,G}{p-q+4u} = \frac{2u\,G - J}{2u-1} ,
\label{G(u)}
}
and, since $p\ne q$,
\equ{
\tr\big(G(u)J\big) = \frac{(p-q)(p+q-4u)}{p-q+4u}
= \frac{4u-2n}{2u-1} = \frac{2\,(2u-\ka-1)}{2u-1} \;\ne\; 0 .
\label{trGJ}
}
This is the first point of departure from the split type DI of \cite{Re26}, where $p=q$ and
$\tr(G(u)J)=0$. Moreover, the term $\ga^{(r)}_{ij}$ of
\eqref{SY-graded-iso} does not vanish here: expanding \eqref{G(u)} gives
$G(u) = G + \tfrac{q-p}{4}(G-J)\,u^{-1} + \mc{O}(u^{-2})$, so that \eqref{S->TGT} together with
$T(u) = I + \bb{F}u^{-1}+\mc{O}(u^{-2})$, where $\bb{F} := \sum_{i,j}E_{ij}\ot F_{ij}$, yields, for the coefficient $S^{(1)}$ of $u^{-1}$
in $S(u)$,
\equ{
S^{(1)} = \bb{F}\,G - G\,\bb{F}^t + \tfrac{q-p}{4}\,(G-J) ,
\label{S1}
}
whose $(i,j)$ entry is $F^{(\rho,0)}_{ij}+\ga^{(0)}_{ij}$. For the pair \eqref{pair} we have
$\tfrac{q-p}{4}=\tfrac12$.


\subsection{The matrix \texorpdfstring{$\wh{S}(u)$}{S(u)} and its Gaussian decomposition} \label{sec:Shat}

The second point of departure is that $S(u)$ admits no Gaussian decomposition
$S(u) = F(u)D(u)E(u)$, with $F(u)$ lower unitriangular, $E(u)$ upper unitriangular and $D(u)$
diagonal with invertible entries, in any ordering of the basis of $\C^N$. Such a decomposition
would make every leading principal submatrix of $G=S(\infty)$ invertible. But $G$ is the
permutation matrix of the transposition $\si=(n,n+1)$. In any ordering, the first leading
principal submatrix that contains $n$ or $n+1$ contains only one of them, so it has a zero row.

Accordingly we work throughout with the adjusted generating matrix
\equ{
\wh{S}(u) := G\,S(u) , \qq \wh{S}(\infty) = G^2 = I ,
\label{Shat}
}
whose leading principal quasi-minors all have leading term $1$ and are therefore invertible. By
\cite[Thm.~4.9.7]{GGRW05}, whose proof uses only this invertibility, it admits a unique Gaussian
decomposition
\equ{
\wh{S}(u) = \wh{F}(u)\,\wh{D}(u)\,\wh{E}(u)
\label{Shat-FDE}
}
over $X^\tw(\mfo_{2n},G)[[u^{-1}]]$, with $\wh{D}(u) = \diag(\wh{d}_1(u),\ldots,\wh{d}_N(u))$,
$\wh{F}(u)$ lower unitriangular with entries $\wh{f}_{ji}(u)$ for $j>i$ and $\wh{E}(u)$
upper unitriangular with entries $\wh{e}_{ij}(u)$ for $i<j$; explicitly, these are given by the
quasi-determinantal formulas of \cite[\S3.1]{Re26} applied to $\wh{S}(u)$ in place of $S(u)$. We
write
\equ{
\wh{f}_{ji}(u) = \sum_{r\ge0}\wh{f}^{(r)}_{ji}u^{-r-1} , \qq
\wh{d}_i(u) = 1 + \sum_{r\ge0}\wh{d}_{i,r}u^{-r-1} , \qq
\wh{e}_{ij}(u) = \sum_{r\ge0}\wh{e}^{(r)}_{ij}u^{-r-1} ,
}
and call these the {\it Gaussian series} and their coefficients the {\it Gaussian generators}.

Since $g_{1a} = \del_{1a}$ and $g_{aN} = \del_{aN}$ (here $n\ge2$, so $1,N\notin\{n,n+1\}$), we
have $\wh{s}_{1k}(u) = s_{1k}(u)$ and $[\wh{S}^{-1}(u)]_{kN} = [S^{-1}(u)]_{kN}$, whence
\equ{
s_{11}(u) = \wh{d}_1(u) , \qq \big[S^{-1}(u)\big]_{NN} = \wh{d}^{\,-1}_N(u) .
\label{corner}
}

\begin{remark} \label{R:not-RE}
The matrix $\wh{S}(u)$ does {\it not}  satisfy the twisted reflection equation \eqref{RE}. Every relation among the Gaussian series below is therefore obtained either by rewriting \cite[(2.19)]{Re26} in terms of $\wh{S}(u)$, or, as we do below,
through the reduction of Section~\ref{sec:emb}, the low-rank case of Section~\ref{sec:low},
the type AI block of Section~\ref{sec:iso}, and the degenerate block of Section~\ref{sec:deg}.
\end{remark}

Here and below $\si=(n,n+1)$ is the transposition above, so $\si(a)=a$ for $a\ne n,n+1$. The
next lemma lifts $\si$ to an automorphism of the algebra. By \eqref{Dr2}--\eqref{Dr3} this
automorphism exchanges the Drinfeld series $b_{n-1}(u)$ and $b_n(u)$ of the two fork nodes. We use
it to obtain the relations of one fork node from those of the other (Sections~\ref{sec:low}
and~\ref{sec:rels}).

\begin{lemma} \label{L:sigma}
The mapping $\vartheta : s_{ab}(u) \mapsto s_{\si(a)\si(b)}(u)$ defines an involutive automorphism
of $X^\tw(\mfo_{2n},\mfo_{n-1}\op\mfo_{n+1})$. For $1\le a\le n-1$ and $b>a$ it satisfies
$\vartheta(\wh{d}_a(u)) = \wh{d}_a(u)$ and $\vartheta(\wh{f}_{ba}(u)) = \wh{f}_{\si(b)a}(u)$.
\end{lemma}

\begin{proof}
Let $P_\si$ be the permutation matrix of $\si$. Since $\ol{\si(a)} = \si(\bar{a})$ for all $a$,
conjugation by $P_\si\ot P_\si$ fixes $P$ and $Q$, and hence $R(u)$ and $R'(u)$. Conjugating
\eqref{RE} by $P_\si\ot P_\si$ shows that the matrix $P_\si\,S(u)\,P_\si$, whose entries are the
series $s_{\si(a)\si(b)}(u)$, satisfies \eqref{RE}. Its constant part is $P_\si\,G\,P_\si = G$.
Hence $\vartheta$ is a homomorphism, and $\vartheta^2=\id$ because $\si^2=\id$.

Next, $\vartheta(\wh{S}(u)) = G\,P_\si\,S(u)\,P_\si = P_\si\,\wh{S}(u)\,P_\si$, that is,
$\vartheta(\wh{s}_{ab}(u)) = \wh{s}_{\si(a)\si(b)}(u)$. By the quasi-determinantal formulas of
\cite[\S3.1]{Re26}, $\wh{d}_a(u)$ and $\wh{f}_{ba}(u)\,\wh{d}_a(u)$ are the quasi-determinants of the
submatrices of $\wh{S}(u)$ on the rows $\{1,\ldots,a\}$, respectively $\{1,\ldots,a-1,b\}$, and
the columns $\{1,\ldots,a\}$, at their last entry. For $a\le n-1$ the permutation $\si$ fixes
$1,\ldots,a$, so $\vartheta$ maps these quasi-determinants to the same ones with $b$ replaced by
$\si(b)$.
\end{proof}


\subsection{Central series and the centre} \label{sec:centre}

Two central series are attached to the symmetry and unitarity relations. We recall their
properties from \cite[\S\S4,\,5]{GR16}. They are proved there from the reflection equation alone and
hold for every $G$ of type DI in \cite[\S2.2]{Re26}. The first is the series
$\ms{c}^\circ(u)\in 1+u^{-1}X^\tw(\mfo_{2n},G)[[u^{-1}]]$ defined by
\[
Q\,S_1(u)\,J_1\,R(2u-\ka)\,J_2\,S^{-1}_2(\ka-u) = \ms{p}^\circ(u)\,\ms{c}^\circ(u)\,Q \qu\text{with}\qu \ms{p}^\circ(u) = 1 - \frac{1}{2u-\ka} + \frac{\tr(G(u)J)}{2u-2\ka} .
\]
Its coefficients are central, it satisfies $\ms{c}^\circ(\ka-u)=\ms{c}^{\circ\,-1}(u)$, and the
symmetry relation \eqref{sym} is equivalent to $\ms{c}^\circ(u)=1$. The second is the series $\ms{z}(u)$ defined by
\equ{
S(u)\,J\,S(-u)\,J = \ms{z}(u)\cdot I .
\label{SJSJ=z}
}
It is a scalar series with central coefficients, and the unitarity
relation \eqref{unit} is the identity $\ms{z}(u)=1$. We need both series in terms of the Gaussian generators of
$\wh{S}(u)$. The first acquires the rational prefactor $\ms{p}(u)$ below because
$\tr(G(u)J)\ne0$ by \eqref{trGJ}. In the split type DI, where $p=q$, this trace is zero and there is no prefactor \cite[(5.16)]{Re26}.

\begin{lemma} \label{L:c(u)}
Set
\equ{
\ms{p}(u) := \frac{(2u-1)(u-\ka)}{u\,(2u-2\ka+1)} , \qq
\ms{c}(u) := \ms{p}(u)\,\wh{d}_1(u)\,\wh{d}^{\,-1}_1(\ka-u) = 1 + \sum_{r\ge0}\ms{c}_r\,u^{-r-1} .
\label{c(u)}
}
Then $\ms{c}(u) = \ms{c}^\circ(u)$ and $\ms{p}(u)\,\ms{p}(\ka-u)=1$. In particular the
coefficients of $\ms{c}(u)$ are central, $\ms{c}(\ka-u)=\ms{c}^{-1}(u)$, and the symmetry
relation \eqref{sym} is equivalent to $\ms{c}(u)=1$.
\end{lemma}

\begin{proof}
By \cite[(5.23)]{GR16}, in the normalisation of \cite[(2.22)]{Re26}, the $(1,1)$ entry of the
 symmetry relation of $X^\tw(\mfo_{2n},G)$ reads
$s_{11}(u)\big(1-\frac{1}{2u-\ka}\big) = \ms{p}^\circ(u)\,\ms{c}^\circ(u)\,s_{11}(\ka-u)$,
since $\del_{1N}=0$. By \eqref{trGJ},
\[
\ms{p}^\circ(u) = \frac{(2u-\ka-1)\,u\,(2u-2\ka+1)}{(2u-\ka)(2u-1)(u-\ka)}
= \frac{1}{\ms{p}(u)}\Big(1-\frac{1}{2u-\ka}\Big) ,
\]
and $s_{11}(u)=\wh{d}_1(u)$ by \eqref{corner}, so $\ms{c}^\circ(u)=\ms{c}(u)$. The identity
$\ms{p}(u)\,\ms{p}(\ka-u)=1$ is a direct computation.
\end{proof}

\begin{lemma} \label{L:z(u)}
The unitarity series is even and is given by
\equ{
\ms{z}(u) = \wh{d}_1(-u)\,\wh{d}_N(u) = \wh{d}_1(u)\,\wh{d}_N(-u) = 1 + \sum_{r\ge0}\ms{z}_r\,u^{-r-1} .
\label{z(u)}
}
\end{lemma}

\begin{proof}
Comparing the $(N,N)$ entries of $J\,S(-u)\,J = \ms{z}(u)\,S^{-1}(u)$ and using \eqref{corner}
gives $\ms{z}(u)=\wh{d}_1(-u)\,\wh{d}_N(u)$. Conjugating \eqref{SJSJ=z} by
$S(u)J$ and using centrality gives $S(-u)\,J\,S(u)\,J = \ms{z}(u)\cdot I$, that is,
$\ms{z}(-u)=\ms{z}(u)$; substituting $u\mapsto-u$ in the expression just obtained yields the
second equality in \eqref{z(u)}.
\end{proof}

The following is refined in Corollary \ref{C:iso-ext} below.

\begin{prop} \label{P:centre}
Let $n\ge2$. The elements $\ms{c}_{2r}$ and $\ms{z}_{2r+1}$, $r\ge0$, are algebraically
independent and generate the centre of $X^\tw(\mfo_{2n},\mfo_{n-1}\op\mfo_{n+1})$, while
$\ms{c}_{2r+1}$ is a polynomial in $\ms{c}_0,\ldots,\ms{c}_{2r}$ and $\ms{z}_{2r}=0$. Moreover
\equ{
X^\tw(\mfo_{2n},\mfo_{n-1}\op\mfo_{n+1}) \;\cong\;
SY^\tw(\mfo_{2n},\mfo_{n-1}\op\mfo_{n+1})\ot\C[\ms{c}_0,\ms{c}_2,\ldots]\ot
\C[\ms{z}_1,\ms{z}_3,\ldots] .
\label{X=SYcz}
}
\end{prop}

\begin{proof}
The polynomial dependence of $\ms{c}_{2r+1}$ follows by comparing coefficients in $\ms{c}(u)\,\ms{c}(\ka-u)=1$.
Lemma~\ref{L:z(u)} yields $\ms{z}_{2r}=0$. The rest is \cite[Cors.~5.3, 5.4]{GR16} for the type BDI matrix of
their \S3.1 with $p-q=2$\footnote{\cite{GR16} uses an opposite convention for $p$ and
$q$: $p\ge q$.}, which is the matrix $G$ of \eqref{G} transported by the substitution of
\cite[\S2.4]{Re26}. In the notation of \cite{GR16}, our $\ms{c}_{2r}$ is $\ms{c}_{2r+1}$, and a
triangular change of central generators replaces $\ms{w}_{2r+2}$ by $\ms{z}_{2r+1}$.
\end{proof}


\subsection{Normalisation invariance} \label{sec:norm}

Some relations below are first proved between the Gaussian images under the homomorphism of
Proposition \ref{P:o4o3->gl2}. Its kernel is generated by the coefficients of $\ms{c}(u)-1$. These
relations therefore hold in the quotient by this kernel. They also hold in its further quotient
$SY^\tw$. To lift them to $X^\tw$, we use the normalised subalgebra of
\cite[Thms.~4.2, 5.3]{GR16}. We follow their \S\S4.2,\,5.1.

Call a relation among the Gaussian series \eqref{Shat-FDE} {\it normalisation invariant} if it is
unchanged when every $\wh{d}_i(u)$ is replaced by $\mu(u)\,\wh{d}_i(u)$, for any series
$\mu(u) = 1 + \sum_{r\ge0} \mu_r\,u^{-r-1}$ with central coefficients $\mu_r\in X^\tw(\mfo_{2n},G)$.
In particular all relations among the Drinfeld series \eqref{Dr1}--\eqref{Dr3} are normalisation
invariant.

Let $\ms{v}(u)$ and $\ms{q}(u)$ be the unique series with constant terms $1$ and central
coefficients such that $\ms{v}^2(u)=\ms{c}(u)$, $\ms{q}(\ka-u)=\ms{q}(u)$ and
$\ms{q}(u)\,\ms{q}(-u)=\ms{v}^{-1}(u)\,\ms{v}^{-1}(-u)\,\ms{z}(u)$. Put
\equ{
\ms{g}(u) := \ms{v}(u)\,\ms{q}(u) .
\label{g(u)}
}
Here $\ms{v}(\ka-u)=\ms{v}^{-1}(u)$ by Lemma \ref{L:c(u)}. The conditions on $\ms{q}(u)$ give
$\ms{q}(u)\,\ms{q}(u+\ka)=\ms{v}^{-1}(u)\,\ms{v}^{-1}(-u)\,\ms{z}(u)$. This equation determines
$\ms{q}(u)$ recursively, and each new coefficient appears with factor $2$. Its right hand side is
even. The equation at $-u$ shows that $\ms{q}(\ka-u)$ solves the same equation. Uniqueness gives
$\ms{q}(\ka-u)=\ms{q}(u)$.
Hence $\ms{g}(u)\,\ms{g}^{-1}(\ka-u)=\ms{c}(u)$ and $\ms{g}(u)\,\ms{g}(-u)=\ms{z}(u)$. So
$\ms{g}^{-1}(u)\,S(u)$ satisfies \eqref{sym} and \eqref{unit}, by Lemma \ref{L:c(u)} and
\eqref{SJSJ=z}.

The subalgebra $\mc{S}$ of $X^\tw(\mfo_{2n},G)$ generated by the
coefficients of the entries of $\ms{g}^{-1}(u)\,S(u)$ is mapped isomorphically onto
$SY^\tw(\mfo_{2n},G)$ by the quotient homomorphism. This is the first tensor factor in
\eqref{X=SYcz}: the multiplication map $\mc{S}\ot\mc{Z}\to X^\tw(\mfo_{2n},G)$ is a linear
isomorphism, where $\mc{Z}$ is the centre.

\begin{lemma} \label{L:X-vs-SY}
The Gaussian series of $G\,\ms{g}^{-1}(u)S(u)$ are $\wh{f}_{ji}(u)$, $\ms{g}^{-1}(u)\,\wh{d}_i(u)$
and $\wh{e}_{ij}(u)$. Consequently a normalisation invariant relation holds in
$X^\tw(\mfo_{2n},G)$ as soon as it holds in $SY^\tw(\mfo_{2n},G)$.
\end{lemma}

\begin{proof}
The first assertion follows by multiplying \eqref{Shat-FDE} by the central scalar
$\ms{g}^{-1}(u)$ and using uniqueness of the Gaussian decomposition. Transfer the relation
through $\mc{S}\cong SY^\tw(\mfo_{2n},G)$ and then remove this scalar by normalisation invariance.
\end{proof}


\subsection{Reduction to lower rank} \label{sec:emb}

We now show that the family \eqref{pair} is closed under the quasi-determinantal reduction of
\cite[\S2.8]{Re26}, that the reduction is compatible with the substitution $S\rightsquigarrow
\wh{S}$, and that it terminates at the pair $(\mfo_4,\mfo_3)$.

Recall from \cite[(2.34)]{Re26} the series
\equ{
s^{[1]}_{ij}(u) := s_{ij}(u) - s_{i1}(u)\,s^{-1}_{11}(u)\,s_{1j}(u) , \qq 1<i,j<N ,
\label{sigma}
}
which are defined because $s_{11}(u) = \wh{d}_1(u)$ is invertible by \eqref{corner}. They are
written $\si_{ij}(u)$ in {\it loc.~cit.}; here $\si$ denotes the transposition $(n,n+1)$.

For $n\ge3$, \cite[Prop.~2.2]{Re26} gives an injective homomorphism
$s_{ij}(u)\mapsto s^{[1]}_{i+1,j+1}(u)$. The remark following that proposition extends it to every
$G$ with $g_{11}=1$. The proof uses the reflection equation, the invertibility of $s_{11}(u)$ and
their Lem.~2.1\,(2), which excludes the target $\mfo_4$. For injectivity it also uses the graded
isomorphism \eqref{X-graded-iso}, $g_{1a}=g_{a1}=\del_{1a}$, and the fact that $G$, $J$ and
$\ga^{(r)}_{ij}$ restrict to those of the source. All of this holds for \eqref{G}, since $q-p=2$ at
every step. The source is again a member of the family. Deleting the rows and
columns $1$ and $N$ from \eqref{G} leaves
\[
G^{[1]} = \sum_{i=2}^{n-1}\big(E_{ii}+E_{\ol{\imath}\,\ol{\imath}}\big) + E_{n,n+1}+E_{n+1,n}
\]
on the span of $e_2,\ldots,e_{N-1}$. After the relabelling $i\mapsto i-1$ this is the matrix
\eqref{G} for the lower rank pair $(\mfo_{2n-2},\mfo_{n-2}\op\mfo_n)$, and
$\ka^{[1]} = \frac{N-2}{2}-\theta = n-2$.

Iterating this reduction by Sylvester's identity, as in \cite[Cor.~2.4]{Re26}, gives
the following.

\begin{prop} \label{P:emb1}
Let $n\ge3$. For $1\le m\le n-2$ the mapping $\psi^{(N)}_m$ of \cite[(2.59)]{Re26}
defines an injective homomorphism
\equ{
\psi^{(N)}_m \;:\; X^\tw\big(\mfo_{2n-2m},\,\mfo_{n-1-m}\op\mfo_{n+1-m}\big) \longrightarrow
X^\tw\big(\mfo_{2n},\,\mfo_{n-1}\op\mfo_{n+1}\big) ,
\label{psi_m}
}
with $\ka^{[m]} = n-1-m$. At $m=n-2$ the source is $X^\tw(\mfo_4,\mfo_3)$ and
$\ka^{[n-2]}=1$. \qed
\end{prop}

For $1\le m\le n-2$ we write $X^{\tw[m]}$ for the image of $\psi^{(N)}_m$ and
$s^{[m]}_{ab}(u) := \psi^{(N)}_m\big(s_{a-m,b-m}(u)\big)$, $m<a,b<\ol{m}$, for the entries of the
reduced generating matrix $S^{[m]}(u) := \sum_{a,b} E_{ab}\ot s^{[m]}_{ab}(u)$.
By \cite[(2.59)]{Re26}, $s^{[m]}_{ab}(u)$ is the
quasi-determinant of the submatrix of $S(u)$ on the rows $\{1,\ldots,m,a\}$ and the columns
$\{1,\ldots,m,b\}$, at the entry $(a,b)$; for $m=1$ it is \eqref{sigma}.
This quasi-determinant is also defined for $m=n-1$,
because the leading $(n-1)\times(n-1)$ submatrix of $S(u)$ has constant term $I$. This gives four
series $s^{[n-1]}_{ab}(u)$ with $a,b\in\{n,n+1\}$. Set
$S^{[n-1]}(u):=\sum_{a,b=n}^{n+1}E_{ab}\ot s^{[n-1]}_{ab}(u)$. We write $X^{\tw[n-1]}$ for the
subalgebra generated by the coefficients of these series. It enters only through
Lemmas \ref{L:Shat-red} and \ref{L:emb-comm}. We need it for the relation
$[\wh{d}_n(u),\wh{f}_{n-1,n-2}(v)]=0$ of Lemma \ref{L:dn-comm}. The series
$\wh{f}_{n-1,n-2}(u)$ has the index $n-2$ outside the band $n-2<a<\ol{n-2}$. So the low rank case
of Section~\ref{sec:low} does not give this relation, and neither does Lemma \ref{L:emb-comm} with
$m=n-2$.

The following lemma makes the substitution $S\rightsquigarrow\wh{S}$ compatible with the
reduction. Let
\equ{
\wh{S}^{[m]}(u) := \wh{F}^{[m]}(u)\wh{D}^{[m]}(u)\wh{E}^{[m]}(u)
}
denote the reduced matrix built, as in \cite[(3.7)--(3.9)]{Re26}, from the Gaussian series of
$\wh{S}(u)$ carrying indices in the band $m<i,j<\ol{m}$.

\begin{lemma} \label{L:Shat-red}
For $1\le m\le n-1$ we have $\wh{S}^{[m]}(u) = G^{[m]}\,S^{[m]}(u)$, where $G^{[m]}$ is the
matrix \eqref{G} of the pair $(\mfo_{2n-2m},\mfo_{n-1-m}\op\mfo_{n+1-m})$; for $m=n-1$ the
matrix $G^{[n-1]}$ is the $2\times2$ antidiagonal permutation matrix. Consequently the
Gaussian series of $\wh{S}(u)$ with indices in the band $m<i,j<\ol{m}$ are the Gaussian series of
the matrix $\wh{S}^{[m]}(u)$ of $X^{\tw[m]}$, and for $m\le n-2$ we have
$X^{\tw[m]} \cong X^\tw(\mfo_{2n-2m},\mfo_{n-1-m}\op\mfo_{n+1-m})$.
\end{lemma}

\begin{proof}
Put $A:=\{1,\ldots,m\}$ and
$B:=\{m+1,\ldots,N-m\}$. We suppress the argument $u$ and write $S_{AB}$ for the submatrix of
$S(u)$ on the rows $A$ and the columns $B$, and similarly for $\wh{S}(u)$. The block $S_{AA}$ has
constant term $I$, so it is invertible. The matrix $G$ fixes $A$ pointwise and preserves $B$, and
$G_{BB}=G^{[m]}$. Hence
\equ{
\wh{S}^{[m]} = \wh{S}_{BB}-\wh{S}_{BA}\,\wh{S}^{-1}_{AA}\,\wh{S}_{AB} = G^{[m]}\big(S_{BB}-S_{BA}\,S^{-1}_{AA}\,S_{AB}\big) = G^{[m]}S^{[m]} .
\label{schur}
}
The first equality holds because the block of $\wh{S}(u)$ on the rows and columns $A\cup B$ is the
product of the corresponding blocks of the factors in \eqref{Shat-FDE}, as in
\cite[Prop.~3.1]{Re26}. The second uses $\wh{S}_{AA}=S_{AA}$, $\wh{S}_{AB}=S_{AB}$,
$\wh{S}_{BA}=G^{[m]}S_{BA}$ and $\wh{S}_{BB}=G^{[m]}S_{BB}$. The last is the quasi-determinantal
formula for $S^{[m]}(u)$ recalled after Proposition \ref{P:emb1}. The remaining assertions follow from
the uniqueness of the Gaussian decomposition and Proposition \ref{P:emb1}.
\end{proof}

\begin{remark} \label{R:no-shift}
Lemma \ref{L:Shat-red} involves no shift of the spectral parameter. Since
$\ka^{[n-2]}=1$ for every $n\ge2$, the Gaussian series of $\wh{S}(u)$ carrying indices in the
band $n-2 < i,j < \ol{n-2}$ satisfy the relations of $X^\tw(\mfo_4,\mfo_3)$ derived in
Section~\ref{sec:low} verbatim, with the index shift $i\mapsto i+n-2$ and with the same spectral
parameter.
\end{remark}


\subsection{The type AI block} \label{sec:iso}

For $1\le i,j,k,l\le n$, restricting \cite[(2.19)]{Re26} gives the relations of type AI:
\eqa{
[\,s_{ij}(u),\,s_{kl}(v)\,] &= \frac{s_{kj}(u)\,s_{il}(v)-s_{kj}(v)\,s_{il}(u)}{u-v}
\\[0.2em]
& \qu - \frac{s_{ik}(u)\,s_{jl}(v)-s_{ki}(v)\,s_{lj}(u)}{u+v-\ka}
+ \frac{s_{ki}(u)\,s_{jl}(v)-s_{ki}(v)\,s_{jl}(u)}{(u-v)(u+v-\ka)} .
\label{RE-iso}
}
That is, the leading $n\times n$ submatrix of $S(u)$ satisfies the twisted reflection equation
\cite[(5.3)]{Re26}. Its constant term is $\diag(1,\ldots,1,0)$, because $g_{nn}=0$. So below we use
only its leading $(n-1)\times(n-1)$ block. The $n\times n$ block returns at $n=3$ in
Section~\ref{sec:deg}.

Since $g_{ia}=\del_{ia}$ for $i<n$, the rows $1,\ldots,n-1$ of $\wh{S}(u)=G\,S(u)$ and of $S(u)$
coincide. Hence the quasi-determinantal formulas of \cite[\S3.1]{Re26} express the Gaussian
series $\wh{d}_i(u)$, $\wh{f}_{ji}(u)$, $\wh{e}_{ij}(u)$ with $1\le i,j\le n-1$ in the entries of
the submatrix $\big(s_{ij}(u)\big)_{1\le i,j\le n-1}$.
Together with \eqref{RE-iso} this gives the following.

\begin{lemma} \label{L:AI-block}
The assignment $s_{ij}(u)\mapsto s_{ij}(u+\ka/2)$, $1\le i,j\le n-1$, defines a homomorphism from
the extended twisted Yangian $X^+(\mfgl_{n-1})$, which is $X_G(\mfo_{n-1})$ with $G=I$ in
\cite[Def.~2.13.1]{Mo07}, onto the subalgebra of
$X^\tw(\mfo_{2n},\mfo_{n-1}\op\mfo_{n+1})$ generated by the coefficients of the series $s_{ij}(u)$
with $1\le i,j\le n-1$. We denote this subalgebra by $X^+_\ka(\mfgl_{n-1})$. Moreover the
Gaussian series $\wh{d}_i(u)$, $\wh{f}_{ji}(u)$, $\wh{e}_{ij}(u)$ of $\wh{S}(u)$ with
$1\le i,j\le n-1$ coincide with the Gaussian series of the generating matrix
$\big(s_{ij}(u)\big)_{1\le i,j\le n-1}$ of $X^+_\ka(\mfgl_{n-1})$. \qed
\end{lemma}

Consequently the whole of \cite[Lem.~5.1, Rem.~5.2]{Re26}, which is derived from the defining
relations of $X^+(\mfgl_n)$ of \cite{LWZ25a}, applies here with $n$ replaced by $n-1$ in the index ranges. The constant $\ka$ remains $n-1$. We
record the statements we use, all for $1\le i,j\le n-2$ and $1\le k,l\le n-1$:
\ali{
[\wh{d}_k(u),\wh{d}_l(v)] = 0 , & \qq [\wh{d}_k(u),\wh{f}_{i+1,i}(v)]=0 \qu\text{if}\qu k\ne i,i+1 ,
\label{AI:dd}
\\
\wh{e}_{i,i+1}(u) = \wh{f}_{i+1,i}(\ka-i-u) , & \qq
\wh{d}_i(u)\,\wh{d}_{i+1}(\ka-i-u) = \wh{d}_{i+1}(u)\,\wh{d}_i(\ka-i-u) ,
\label{AI:f=e}
}
together with \cite[(5.5)--(5.10)]{Re26}, which we do not repeat. The coefficient of $u^{-2}$ in
their (5.5) gives
\equ{
[\,\wh{d}_{i,1},\,\wh{f}_{i+1,i}(v)\,] = (\ka-i-2v)\,\wh{f}_{i+1,i}(v) + 2\,\wh{f}^{\,(0)}_{i+1,i} .
\label{AI:df}
}

We now relate this block to the reduced subalgebras. Let $1\le m\le n-1$. By the
quasi-determinantal formulas noted before Lemma \ref{L:AI-block}, the Gaussian series of
$\wh{S}(u)$ with indices at most $m$ are expressed in the series $s_{ab}(u)$ with $1\le a,b\le m$.
By \cite[Lem.~2.6]{Re26}, the coefficients of these series $s_{ab}(u)$ commute with $X^{\tw[m]}$.
That lemma is stated for split pairs and for all $1\le m<n$, in terms of the quasi-determinants
$s^{[m]}_{ab}(u)$. Its proof is an induction on $m$ through $\psi^{(N)}_1$. It uses $G$ only through
its rows and columns $1,\ldots,m$, which for $m\le n-1$ are those of the identity matrix. At $m=n-1$
the induction ends with their Lem.~2.1\,(1) in $X^\tw(\mfo_4,\mfo_3)$, and this identity holds at
$\ka=1$. This gives the following.

\begin{lemma} \label{L:emb-comm}
Let $1\le m\le n-1$. The coefficients of the series $\wh{d}_a(u)$ with $1\le a\le m$ and
$\wh{f}_{a+1,a}(u)$ with $1\le a<m$ commute with $X^{\tw[m]}$. \qed
\end{lemma}


\section{The low rank case $X^\tw(\mfo_4,\mfo_3)$} \label{sec:low}

By Proposition \ref{P:emb1} the reduction of Section~\ref{sec:emb} ends at $X^\tw(\mfo_4,\mfo_3)$, the
case $n=2$ of \eqref{pair}. Its relations, derived in this section, are transferred to general $n$
in Sections~\ref{sec:deg} and~\ref{sec:Dr} by Lemma \ref{L:Shat-red}.


\subsection{The Yangian $Y(\mfgl_2)$}

Consider the Yangian $Y(\mfgl_2)$ and write its generating matrix $T(u)$ as
\equ{
T(u) = \ms{F}(u)\,\ms{D}(u)\,\ms{E}(u) =
\begin{bmatrix} 1 & 0 \\ \ms{f}(u) & 1 \end{bmatrix}
\begin{bmatrix} \ms{d}_1(u) & 0 \\ 0 & \ms{d}_2(u) \end{bmatrix}
\begin{bmatrix} 1 & \ms{e}(u) \\ 0 & 1 \end{bmatrix} .
\label{Y:FDE}
}
Defining relations of $Y(\mfgl_2)$ then imply the Ding--Frenkel relations \cite[\S3.1]{Mo07}:
\eqg{
\relax[\ms{d}_i(u),\,\ms{d}_j(v)] = 0 , \qq [\ms{e}(u),\,\ms{f}(v)] = \frac{\xi(u)-\xi(v)}{u-v} ,
\\
[\ms{d}_i(u),\,\ms{e}(v)] = (-1)^i\,\frac{\ms{d}_i(u)\,(\ms{e}(u)-\ms{e}(v))}{u-v} , \qq [\ms{e}(u),\,\ms{e}(v)] = \frac{(\ms{e}(u)-\ms{e}(v))^2}{u-v} ,
\\
[\ms{d}_i(u),\,\ms{f}(v)] = (-1)^i\,\frac{(\ms{f}(u)-\ms{f}(v))\,\ms{d}_i(u)}{v-u} , \qq [\ms{f}(u),\,\ms{f}(v)] = \frac{(\ms{f}(u)-\ms{f}(v))^2}{v-u} ,
\label{Y:gl2}
}
where $\xi(u) = \ms{d}_1^{-1}(u)\,\ms{d}_2(u)$ and $i,j \in \{ 1,2\}$.
We shall also use the {\it quantum determinant}
\equ{
\operatorname{qdet}T(u) := \ms{d}_1(u)\,\ms{d}_2(u-1) .
\label{Y:qdet}
}
The coefficients of $\ms{e}(u)$, $\ms{f}(u)$ and $\xi(u)$ generate the {\it special Yangian}
$Y(\mfsl_2)\subset Y(\mfgl_2)$. We write $x^+(u) := \ms{e}(u)$ and $x^-(u) := \ms{f}(u)$ for its
Drinfeld currents, and $\xi(u)$ is its Cartan current. Relations \eqref{Y:gl2} imply
\eqg{
\relax[\xi(u),\,\xi(v)] = 0 , \qq [\xi(u),\,x^\pm(v)] = \pm \,\frac{\{ \xi(u),\, x^\pm(u)-x^\pm(v)\}}{u-v} ,
\\
[x^\pm(u),\,x^\pm(v)] = \pm\,\frac{(x^\pm(u)-x^\pm(v))^2}{u-v} , \qq
[x^+(u),\,x^-(v)] = \frac{\xi(u)-\xi(v)}{u-v} .
\label{Y:drinfeld}
}
Here $\{ a,b\} = ab+ba$. These relations are defining relations of $Y(\mfsl_2)$ in its
Drinfeld presentation \cite{Dr87}; see \cite[Cor.~3.1.7, Rem.~3.1.8]{Mo07} at $N=2$. Our
$x^\pm(u)$ and $\xi(u)$ are the series $\xi^\mp_1(u)$ and $\kappa_1(u)$ there.


\subsection{The extended twisted Yangian $X^\tw(\mfo_4,\mfo_3)$} \label{sec:o4o3}

We specialise Section~\ref{sec:setup} to $n=2$: $N=4$, $(p,q)=(1,3)$,
\equ{
\theta = 1 , \qu \ka = 1 , \qu
G = E_{11} + E_{23} + E_{32} + E_{44} , \qu J = E_{14} + E_{23} + E_{32} + E_{41} ,
\label{o4o3-data}
}
so that the symmetry prefactor of \eqref{c(u)} is $\ms{p}(u) = \frac{u-1}{u}$ and
\equ{
\ms{c}(u) = \frac{u-1}{u}\,\wh{d}_1(u)\,\wh{d}^{\,-1}_1(1-u) , \qq
\ms{z}(u) = \wh{d}_1(-u)\,\wh{d}_4(u) = \wh{d}_1(u)\,\wh{d}_4(-u) ,
\label{o4o3-cz}
}
by \eqref{c(u)} and \eqref{z(u)}. The special twisted Yangian $SY^\tw(\mfo_4,\mfo_3)$ is the
quotient of $X^\tw(\mfo_4,\mfo_3)$ by the ideal generated by the coefficients of $\ms{c}(u)-1$ and
$\ms{z}(u)-1$, by Lemma \ref{L:c(u)} and \eqref{SJSJ=z}.

Let $V = \C^2 \ot \C^2$ and let (as in \cite[\S4.3]{AMR06})
\equ{
v_1 := e_1 \ot e_1, \qu
v_2 := e_1 \ot e_2, \qu
v_3 := e_2 \ot e_1, \qu
v_4 := -e_2 \ot e_2
\label{v-basis}
}
be a basis of $V$, and let $\equiv$ denote equality of linear maps upon the identification of
$e_i \in \C^4$ with $v_i \in V$. Write $P^\circ$ for the flip of the two tensor factors of $V$ and
$Q^\circ := \sum_{i,j=1}^2 E_{ij} \ot E_{ij}$, and set $R^\circ(u) = I - u^{-1} P^\circ$, so that
$R^{\circ\prime}(u) = I - u^{-1} Q^\circ$, the prime denoting the partial transpose in the first
tensor factor. Then
\equ{
G \equiv P^\circ , \qq Q^\circ \equiv G - J , \qq
G(u) \equiv G\,R^{\circ\prime}(\ka-2u) = G - \frac{Q^\circ}{1-2u} ,
\label{G-flip}
}
in agreement with \eqref{G(u)}. Let $(X,Y)\in\mfsl_2\op\mfsl_2\cong\mfo_4$ act on $V$ by
$X\ot1+1\ot Y$. Then conjugation by $G$ is the flip $(X,Y)\mapsto(Y,X)$, and
$\rho(X,Y)=(-Y^t,-X^t)$. The automorphism $(X,Y)\mapsto(X,-Y^t)$ conjugates $\rho$ to the flip. So
$(\mfo_4,\mfo_3) \cong (\mfsl_2 \op \mfsl_2, \Delta \mfsl_2)$ is the {\it diagonal} symmetric pair.

We now realise the quotient $X^\tw(\mfo_4,\mfo_3)/(\ms{c}(u)-1)$ inside
$Y(\mfgl_2)\ot Y(\mfgl_2)$. This quotient contains a copy of $SY^\tw(\mfo_4,\mfo_3)$ by
Section~\ref{sec:norm}. We write the image of $S(u)$ through one matrix $\wt{B}(u)$ that satisfies
the relations of $Y(\mfgl_2)$, so that relations \eqref{Y:gl2} apply to its Gaussian series.

Consider the tensor product $Y(\mfgl_2) \ot Y(\mfgl_2)$, and let $T(u)$ (resp.~$\ol{T}(u)$)
denote the generating matrix of the first (resp.~second) $Y(\mfgl_2)$ factor. Set
\equ{
\wt B(u) := T(u)\,\ol{T}{}^{t}(-u) .
\label{B(u)}
}
The matrix $\wt{B}(u)$ satisfies the defining relations of $Y(\mfgl_2)$: indeed
$\nu: T(u) \mapsto T^{\,t}(-u)$ is an automorphism of $Y(\mfgl_2)$, and $\wt{B}(u)$ is the
image of $T(u)$ under $(\id\ot\nu)\circ\Delta$. Since
$(\id\ot\eps)\circ(\id\ot\nu)\circ\Delta=\id$, this homomorphism is injective, so the
subalgebra $\langle\wt B\rangle$ of $Y(\mfgl_2)\ot Y(\mfgl_2)$ generated by the coefficients of
the entries of $\wt B(u)$ is a copy of $Y(\mfgl_2)$. We write
\equ{
\wt{B}(u) = \wt{\ms{F}}(u)\,\wt{\ms{D}}(u)\,\wt{\ms{E}}(u)
\label{B-gauss}
}
for its Gaussian decomposition \eqref{Y:FDE} and the associated series
$\wt{\ms{d}}_1(u),\wt{\ms{d}}_2(u),\wt{\ms{e}}(u),\wt{\ms{f}}(u)$ and
$\wt\xi(u) := \wt{\ms{d}}{}^{-1}_1(u)\,\wt{\ms{d}}_2(u)$, which obey exactly the same relations
\eqref{Y:gl2}.

Next, we set
\equ{
\wt{S}_{12}(u) := \wt{B}_2(u-\tfrac12)\, R^{\circ\prime}_{12}(1-2u)\,\wt{B}^t_1(\tfrac12-u) ,
\label{S(u)-diag}
}
the subscripts referring to the two $\C^2$ factors of the auxiliary space $V$. Insert
\eqref{B(u)}. The two Yangian factors commute, so
$\wt{B}^t_1(\tfrac12-u)=\ol T_1(u-\tfrac12)\,T^t_1(\tfrac12-u)$. Partial transposition of the $RTT$
relation of $Y(\mfgl_2)$ gives
\[
T_1(a)\,R^{\circ\prime}_{12}(b-a)\,T^t_2(b) = T^t_2(b)\,R^{\circ\prime}_{12}(b-a)\,T_1(a) .
\]
Apply this identity to $\ol T(u)$ with $a=u-\tfrac12$ and $b=\tfrac12-u$. Then $b-a=1-2u$. The
identity interchanges the middle three of the five factors of $\wt{S}_{12}(u)$. The last two factors
commute because they lie in different Yangian and auxiliary factors. We obtain
\equ{
\wt{S}_{12}(u) = T_2(u-\tfrac12)\,\ol T_1(u-\tfrac12)\,R^{\circ\prime}_{12}(1-2u)\,T^t_1(\tfrac12-u)\,\ol T^t_2(\tfrac12-u) .
\label{S-two-yangian}
}
\enlargethispage{2em}%
Multiply this by $G$ on the left. Since $G\equiv P^\circ$, this interchanges the subscripts $1,2$
of the first two factors, and the middle factor becomes $G(u)$ by \eqref{G-flip}. Thus
$S(u)\mapsto G\,\wt{S}_{12}(u)$ is the embedding \eqref{S->TGT},
$S(u)\mapsto\bb{T}(u-\tfrac12)\,G(u)\,\bb{T}^{\,t}(\tfrac12-u)$, for the matrix
$\bb{T}(u) := T_1(u)\,\ol T_2(u)$. By \cite[Prop.~4.8]{AMR06}, $\bb{T}(u)$ is the image of the
generating matrix of the extended Yangian $X(\mfo_4)$ under the embedding
$X(\mfo_4)\hookrightarrow Y(\mfgl_2)\ot Y(\mfgl_2)$. So by \cite[Lems.~4.2, 4.3]{GR16}, in the
presentation of \cite[\S2.4]{Re26}, the matrix $G\,\wt{S}(u)$ satisfies the reflection equation
\eqref{RE} and the symmetry relation \eqref{sym}. It does not satisfy the unitarity relation \eqref{unit},
because $\bb{T}(u)$ is that of $X(\mfo_4)$ and not of $Y(\mfo_4)$.

Write $\langle\wt{S}\rangle$ for the subalgebra of $Y(\mfgl_2)\ot Y(\mfgl_2)$ generated by the
coefficients of the entries of $G\,\wt{S}(u)$, and write $\wt{\ms{c}}(u)$ and $\wt{\ms{z}}(u)$
for the series \eqref{c(u)} and \eqref{z(u)} formed from $G\,\wt{S}(u)$. By \eqref{sym} and Lemma \ref{L:c(u)}
we have $\wt{\ms{c}}(u)=1$, while $\wt{\ms{z}}(u)\ne1$ by \eqref{cz-image} below.

\begin{prop}\label{P:o4o3->gl2}
The mapping $X^\tw(\mfo_4,\mfo_3) \to \langle\wt{S}\rangle$ given by
$S(u) \mapsto G\,\wt{S}(u)$, subject to the identification of basis vectors described above, is a
surjective homomorphism of algebras whose kernel is the two-sided ideal generated by the
coefficients of $\ms{c}(u)-1$. Its restriction to the subalgebra $\mc{S}$ of
Section~\ref{sec:norm}, which the quotient homomorphism identifies with $SY^\tw(\mfo_4,\mfo_3)$, is
injective. The Gaussian images needed below are (here $x=u - \tfrac12$)
\eqg{
\wh{d}_1(u) \mapsto \frac{2u}{2u-1}\,\wt{\ms{d}}_1(-x)\,\wt{\ms{d}}_1(x) , \qq
\wh{d}_2(u) \mapsto \wt{\ms{d}}_1(-x)\,\wt{\ms{d}}_2(x) ,
\\
\wh{d}_3(u) \mapsto \wt{\ms{d}}_1(-x)\,\wt\xi(-x-1)\,\wt{\ms{d}}_1(x) , \qq
\wh{d}_4(u) \mapsto \frac{2u+1}{2u}\,\wt{\ms{d}}_1(-x)\,\wt\xi(-x-1)\,\wt{\ms{d}}_2(x) ,
\\
\wh{f}_{21}(u) \mapsto \wt{\ms{f}}(x) , \qq
\wh{f}_{31}(u) \mapsto \wt{\ms{e}}(-x-1) , \qq
\wh{f}_{32}(u) \mapsto 0 ,
\\
\wh{e}_{12}(u) \mapsto \wt{\ms{e}}(x) , \qq
\wh{e}_{13}(u) \mapsto \wt{\ms{f}}(-x-1) , \qq
\wh{e}_{23}(u) \mapsto 0 .
\label{o4o3->gl2}
}
\end{prop}

\begin{proof}
The homomorphism property was shown above, and the mapping is surjective by the definition of
$\langle\wt{S}\rangle$. By \eqref{S-two-yangian} it is the composite of the quotient map
$X^\tw(\mfo_4,\mfo_3)\to X^\tw(\mfo_4,\mfo_3)/(\ms{c}(u)-1)$, of the embedding \eqref{S->TGT} of
the quotient into the extended Yangian $X(\mfo_4)$, and of the embedding
$X(\mfo_4)\hookrightarrow Y(\mfgl_2)\ot Y(\mfgl_2)$ of \cite[Prop.~4.8]{AMR06}. So its kernel is
the ideal generated by the coefficients of $\ms{c}(u)-1$. This ideal lies in the kernel of the
quotient homomorphism $X^\tw(\mfo_4,\mfo_3)\to SY^\tw(\mfo_4,\mfo_3)$, which is injective on
$\mc{S}$ by Section~\ref{sec:norm}. So the mapping is injective on $\mc{S}$.

The images \eqref{o4o3->gl2} of the Gaussian series follow by a direct computation:
Gauss-decompose $\wt S(u)$, the image of $\wh S(u)$, as in \eqref{Shat-FDE}, and simplify with
\eqref{Y:gl2}. The $RTT$ relation at neighbouring arguments gives
\[
\wt{B}^t(y) =
\begin{bmatrix} 1 & 0 \\ \wt{\ms{e}}(y-1) & 1 \end{bmatrix}
\begin{bmatrix} \wt{\ms{d}}_1(y) & 0 \\ 0 & \wt{\ms{d}}_1(y)\,\wt\xi(y-1) \end{bmatrix}
\begin{bmatrix} 1 & \wt{\ms{f}}(y-1) \\ 0 & 1 \end{bmatrix}
\]
Set $y=-x$. The argument $y-1$ in this decomposition gives the reflected argument $-x-1$ in
\eqref{o4o3->gl2}. For the central series we have
\equ{
\ms{c}(u) \mapsto \wt{\ms{c}}(u) = 1 , \qq
\ms{z}(u) \mapsto \wt{\ms{z}}(u)
= \operatorname{qdet}\wt B(u+\tfrac12) \operatorname{qdet}\wt B(\tfrac12-u) .
\label{cz-image}
}
The first relation was noted above. For the second, take the image
of $\wh{d}_1(u)\,\wh{d}_4(-u)$ under \eqref{o4o3->gl2} and apply
$\wt{\ms{d}}_1(x)\,\wt\xi(x) = \wt{\ms{d}}_2(x)$ and
$\operatorname{qdet}\wt B(y) = \wt{\ms{d}}_1(y)\,\wt{\ms{d}}_2(y-1)$. The mapping therefore does
not factor through $SY^\tw(\mfo_4,\mfo_3)$.
\end{proof}

Proposition \ref{P:o4o3->gl2} and the Ding--Frenkel relations \eqref{Y:gl2} give the
wanted identities between the Gaussian series of $\wh{S}(u)$. We write them in terms of the series
\equ{
b_1(u) := \wh{f}_{21}(u) , \qq b_2(u) := \wh{f}_{31}(u) , \qq
h_1(u) := \frac{2u}{2u-1}\,\wh{d}^{\,-1}_1(u)\,\wh{d}_2(u) , \qq h_2(u) := h_1(-u) ,
\label{o4o3:Dr}
}
which are the wanted Drinfeld series \eqref{Dr2}--\eqref{Dr3} at $n=2$.

\begin{crl} \label{C:o4o3}  
The following relations hold in $X^\tw(\mfo_4,\mfo_3)$ with $1\le k,l \le 4$ and $1\le i,j\le 2$:
\gat{
\wh{e}_{12}(u) = b_2(-u) , \qq \wh{e}_{13}(u) = b_1(-u) , \qq \wh{f}_{32}(u) = \wh{e}_{23}(u) = 0 ,
\label{fork:ef+triv}
\\
[\wh{d}_k(u),\wh{d}_l(v)] = 0  , \qq [h_1(u),h_1(v)] = 0 , \qq h_2(u) = \frac{2u}{2u-1}\,\wh{d}^{\,-1}_1(u)\,\wh{d}_3(u) ,
\label{fork:hh}
\\
\wh{d}_1(u)\,b_j(v) = \bigg( \frac{(u-v-1)(u+v)}{(u-v)(u+v-1)}\,b_j(v)
+ \frac{2u}{(2u-1)(u-v)}\,b_j(u)
- \frac{2(u-1)}{(2u-1)(u+v-1)}\,b_j(1-u) \bigg)\,\wh{d}_1(u) ,
\label{d-b}
\\[0.2em]
[h_1(u),b_1(v)] = -\,\frac{\{h_1(u),\,b_1(u)-b_1(v)\}}{u-v} , \qq
[h_1(u),b_2(v)] = \frac{\{h_1(u),\,b_2(-u)-b_2(v)\}}{u+v} ,
\label{fork:hb}
\\
[b_i(u),b_i(v)] = -\,\frac{(b_i(u)-b_i(v))^2}{u-v} , \qq
[b_1(u),b_2(v)] = \frac{h_2(v)-h_1(u)}{u+v} .
\label{fork:bb}
}
\end{crl}

\begin{proof}
By Proposition \ref{P:o4o3->gl2}, $\langle\wt{S}\rangle$ is the quotient of
$X^\tw(\mfo_4,\mfo_3)$ by the coefficients of $\ms{c}(u)-1$. So a relation between the images in
$\langle\wt{S}\rangle$ holds in this quotient, and hence in its further quotient
$SY^\tw(\mfo_4,\mfo_3)$. All the relations above are normalisation invariant. So Lemma
\ref{L:X-vs-SY} lifts them to $X^\tw(\mfo_4,\mfo_3)$.

All relations except \eqref{d-b} follow directly from \eqref{o4o3->gl2} and
\eqref{Y:gl2}. By Lemma \ref{L:sigma} it suffices to prove \eqref{d-b} for $j=1$. By
\eqref{Y:gl2} we have
\[
\wt{\ms{d}}_1(a)\,\wt{\ms{f}}(b) = \bigg( \frac{a-b-1}{a-b}\,\wt{\ms{f}}(b) + \frac{1}{a-b}\,\wt{\ms{f}}(a) \bigg)\,\wt{\ms{d}}_1(a) .
\]
Write $x = u-\tfrac12$ and $y = v-\tfrac12$. By \eqref{o4o3->gl2}, apply this rule to move first
$\wt{\ms{d}}_1(x)$ and then $\wt{\ms{d}}_1(-x)$ past $\wt{\ms{f}}(y)$. The resulting series are the
images of $b_1(v)$, $b_1(u)$ and $b_1(1-u)$, with the coefficients in \eqref{d-b}.
\end{proof}

We now identify $SY^\tw(\mfo_4,\mfo_3)$ with the Yangian $Y(\mfsl_2)$. For this we need the symbols
of the coefficients of $b_1(u)$, $b_2(u)$ and $h_1(u)$ in the associated graded algebra. Write
$\mfo^{(k)}x^k$, with $\mfo^{(k)}\subseteq\mfo_4$, for the degree-$k$ component of
$\mfo_4[x]^\rho$, so that $\gr SY^\tw(\mfo_4,\mfo_3)\cong U(\mfo_4[x]^\rho)$ by
\eqref{SY-graded-iso}. Write also $h_1(u) = 1+\sum_{k\ge0}h_{1,k}\,u^{-k-1}$.

\begin{lemma} \label{L:o4o3-symbols}
For every $k\ge0$ the symbols of $\wh{f}^{\,(k)}_{21}$, $\wh{f}^{\,(k)}_{31}$ and $h_{1,k}$ form a
basis of $\mfo^{(k)}x^k$.
\end{lemma}

\begin{proof}
We compute the three symbols by \eqref{SY-graded-iso} and then count dimensions. Since
$\wh{s}_{ij}(u)=s_{\si(i)j}(u)$, we have
\[
\wh{f}_{21}(u) = s_{31}(u)\,s^{-1}_{11}(u) , \qu
\wh{f}_{31}(u) = s_{21}(u)\,s^{-1}_{11}(u) , \qu
h_1(u) = \frac{2u}{2u-1}\,s^{-1}_{11}(u)\,\big(s_{32}(u)-s_{31}(u)\,s^{-1}_{11}(u)\,s_{12}(u)\big) .
\]
Expand these series in $u^{-1}$. A product of at least two coefficients $s^{(r)}_{ab}$ with $r\ge1$
that contributes to the coefficient of $u^{-k-1}$ has degree at most $k-1$. Also
$\frac{2u}{2u-1}=1+\tfrac12u^{-1}+\mc{O}(u^{-2})$. So the three symbols are
$\bar{s}^{(k+1)}_{31}$, $\bar{s}^{(k+1)}_{21}$ and
$\bar{s}^{(k+1)}_{32}-\bar{s}^{(k+1)}_{11}+\tfrac12\del_{k0}$.

The scalars of \eqref{SY-graded-iso} satisfy $\ga^{(0)}_{31}=\ga^{(0)}_{21}=\ga^{(0)}_{32}=0$
and $\ga^{(0)}_{11}=\tfrac12$. So by \eqref{SY-graded-iso} the three symbols are
$F^{(\rho,k)}_{31}$, $F^{(\rho,k)}_{21}$ and $F^{(\rho,k)}_{32}-F^{(\rho,k)}_{11}$. For the
matrix \eqref{G} with $n=2$ these are
\[
\big(F_{31}-(-1)^kF_{12}\big)x^k , \qq
\big(F_{21}-(-1)^kF_{13}\big)x^k , \qq
\big((1+(-1)^k)F_{33}-(1-(-1)^k)F_{11}\big)x^k ,
\]
which are linearly independent. By \S\ref{sec:notation}, $\mfo^{(k)}$ is the $(-1)^k$-eigenspace of
$\rho$ in $\mfo_4$. Since $\dim\mfo_4=6$ and $\dim\mfo_4^\rho=\dim\mfo_3=3$, we have
$\dim\mfo^{(k)}=3$. So the three elements form a basis of $\mfo^{(k)}x^k$.
\end{proof}

\begin{thm} \label{T:o4o3}
The mapping $x^\pm(u)\mapsto b^\pm(u)$ and $\xi(u)\mapsto h_1(u)$, where $b^-(u) := b_1(u)$ and
$b^+(u) := b_2(-u)$, defines an isomorphism $Y(\mfsl_2) \iso SY^\tw(\mfo_4,\mfo_3)$.
\end{thm}

\begin{proof}
We obtain the homomorphism from Proposition \ref{P:o4o3->gl2} and then show that it is bijective.
By Lemma \ref{L:X-vs-SY} the series $b^\pm(u)$ and $h_1(u)$ are also formed from the Gaussian
series of $G\,\ms{g}^{-1}(u)S(u)$, so their coefficients lie in $\mc{S}$. By Proposition
\ref{P:o4o3->gl2} the mapping \eqref{o4o3->gl2} is injective on
$\mc{S}\cong SY^\tw(\mfo_4,\mfo_3)$ and sends $b^-(u)$, $b^+(u)$, $h_1(u)$ to
$\wt{\ms{f}}(u-\tfrac12)$, $\wt{\ms{e}}(u-\tfrac12)$, $\wt\xi(u-\tfrac12)$. A shift of the argument
does not change \eqref{Y:drinfeld}. Thus $b^-(u)$, $b^+(u)$ and $h_1(u)$ satisfy \eqref{Y:drinfeld}
with $x^\pm(u)=b^\pm(u)$ and $\xi(u)=h_1(u)$. These are defining relations of $Y(\mfsl_2)$, so we
obtain a homomorphism $Y(\mfsl_2)\to SY^\tw(\mfo_4,\mfo_3)$.

Compose it with the inverse of the isomorphism $\mc{S}\cong SY^\tw(\mfo_4,\mfo_3)$ and with the
mapping \eqref{o4o3->gl2} on $\mc{S}$. This composite is the shift of the argument by $-\tfrac12$,
followed by the inclusion of $Y(\mfsl_2)$ into the copy $\langle\wt B\rangle$ of $Y(\mfgl_2)$. So
the homomorphism is injective.

By Lemma \ref{L:o4o3-symbols} the symbols of the coefficients of $b^-(u)$, $b^+(u)$ and $h_1(u)$
span $\mfo^{(k)}x^k$ for every $k$. So they generate
$U(\mfo_4[x]^\rho)\cong\gr SY^\tw(\mfo_4,\mfo_3)$. An induction on the filtration degree shows that
the coefficients generate $SY^\tw(\mfo_4,\mfo_3)$. So the homomorphism is surjective.
\end{proof}

\begin{remark} \label{R:o4o3-vs-split}
Corollary \ref{C:o4o3} uses $b_1(u)$, $b_2(u)$, $h_1(u)$ and $h_2(u)=h_1(-u)$ instead of the
currents of Theorem \ref{T:o4o3}. This makes Theorem \ref{T:iso} uniform in the nodes. The comparison
with the split case is given after that theorem.
\end{remark}


\section{The degenerate block of type AI} \label{sec:deg}

Lemma \ref{L:AI-block} gives the relations between the Gaussian series of $\wh{S}(u)$ with
indices $1\le i,j\le n-1$. Corollary \ref{C:o4o3} does the same for the indices
$n-1\le i,j\le n+2$, through Lemma \ref{L:Shat-red} and Proposition \ref{P:emb1}. The two ranges
overlap in the single index $n-1$, so the diagonal series $\wh{d}_{n-1}(u)$ belongs to both;
together with Lemma \ref{L:emb-comm} the two results give the commutativity of the diagonal
series $\wh{d}_i(u)$ with $1\le i\le n$ (Lemma \ref{L:dd}). They do not give the relations
between $\wh{f}_{n-1,n-2}(u)$ and the series $\wh{f}_{n,n-1}(u)$, $\wh{f}_{n+1,n-1}(u)$ of the second
range.

Let $n\ge3$. We introduce the shorthand notation
\equ{
f(u) := \wh{f}_{n-1,n-2}(u) , \qq b(u) := \wh{f}_{n+1,n-1}(u) , \qq g(u) := \wh{f}_{n+1,n-2}(u) .
\label{deg:notation}
}
We write $f^{(0)}$, $b^{(0)}$ and $g^{(0)}$ for the coefficients of $u^{-1}$ in
$f(u)$, $b(u)$ and $g(u)$.
The series $f(u)$ belongs to the first range and $b(u)$ to the second, while the {\it corner}
series $g(u)$ belongs to neither. We first derive the relations involving $g(u)$, and from them
the relations between $f(u)$ and $b(u)$. The automorphism $\vartheta$ of Lemma \ref{L:sigma}
fixes $f(u)$ and maps $b(u)$ to $\wh{f}_{n,n-1}(u)$. The relations between $f(u)$ and $b(u)$,
together with their images under $\vartheta$, are the last relations needed for Theorem
\ref{T:iso}. %
By the following lemma it suffices to consider the algebra $X^\tw(\mfo_6,\mfo_2\op\mfo_4)$.

\begin{lemma} \label{L:reduce-to-3}
Let $n\ge3$. A relation among the Gaussian series of $\wh{S}(u)$ with indices in the band
$n-3<a<\ol{n-3}$, such as the series $\wh{d}_{i}(u)$ with $n-2\le i \le n+1$ and the series
\eqref{deg:notation}, holds in $X^\tw(\mfo_{2n},\mfo_{n-1}\op\mfo_{n+1})$ if and only if the relation obtained by lowering every index by $n-3$, with the same spectral
parameters and the same coefficients, holds in $X^\tw(\mfo_6,\mfo_2\op\mfo_4)$.
\end{lemma}

\begin{proof}
For $n=3$ there is nothing to prove. For $n\ge4$ apply Lemma \ref{L:Shat-red} with $m=n-3$: these series are the corresponding Gaussian series of
$X^{\tw[n-3]}\cong X^\tw(\mfo_6,\mfo_2\op\mfo_4)$, with $\ka^{[n-3]}=2$ and with no shift of the
spectral parameter. Injectivity of $\psi^{(N)}_{n-3}$ gives the converse.
\end{proof}

Throughout the rest of this section $n=3$, so that $N=6$, $\ka=2$, $G$ is the permutation
matrix of $\si=(3,4)$, and $f(u) = \wh{f}_{21}(u)$, $b(u) = \wh{f}_{42}(u)$, $g(u) = \wh{f}_{41}(u)$.


\subsection{The block} \label{sec:blocks}

Consider the leading $3\times3$ submatrix
\equ{
S^\circ(u) := \sum_{k,l=1}^3 E_{kl}\ot s_{kl}(u)
\label{deg:block}
}
of $S(u)$. By \eqref{RE-iso} it satisfies the twisted reflection equation of type AI with
$\ka=2$. The block is {\it degenerate}: since $g_{33}=0$ we
have $S^\circ(\infty) = \diag(1,1,0)$. Hence $S^\circ(u)$ is not the image of the generating
matrix of $X^+(\mfgl_3)$, whose constant term is the identity matrix. Therefore \cite[Lem.~5.1]{Re26},
which is derived from the relations of $X^+(\mfgl_2)$ and $X^+(\mfgl_3)$ obtained in \cite{LWZ25a},
applies only to the positions $1,2$. They carry the subalgebra $X^+_\ka(\mfgl_2)$ of Lemma \ref{L:AI-block}.

The leading submatrices of $S^\circ(u)$ of sizes $1$ and $2$ have constant terms $1$ and $I$. So the
quasi-determinantal formulas of \cite[\S3.1]{Re26} define the series $d_k(u)$ ($1\le k \le 3$) and
$f_{lk}(u)$, $e_{kl}(u)$ ($1\le k<l\le 3$) of $S^\circ(u)$. No inverse of $d_3(u)$ is needed. We omit
the hats that mark those of $\wh{S}(u)$. The next lemma expresses the series that we need through
the Gaussian series of $\wh{S}(u)$. It also shows that $d_3(u)=0$.

\begin{lemma} \label{L:block}
The series $d_1(u)$ and $d_2(u)$ are invertible, and
\gat{
d_k(u) = \wh{d}_k(u) \qu (k=1,2) , \qq f_{21}(u) = f(u) , \qq e_{12}(u) = \wh{e}_{12}(u) ,
\label{blocks:1}
\\[0.2em]
f_{32}(u) = b(u) , \qq f_{31}(u) = g(u) , \qq e_{23}(u) = b(-u) , \qq d_3(u) = 0 .
\label{blocks:2}
}
\end{lemma}

\begin{proof}
The rows $1,2$ of $\wh{S}(u)$ and of $S(u)$ coincide, as noted before Lemma \ref{L:AI-block}.
This gives \eqref{blocks:1}, and the leading terms of $d_1(u)$ and $d_2(u)$ equal $1$ because
$g_{kl}=\del_{kl}$ for $k,l\le2$.

By \eqref{Shat} the row $3$ of $S(u)$ is the row $4$ of $\wh{S}(u)$. We use the
quasi-determinantal formulas of \cite[\S3.1]{Re26}: for index sets $A$ and $B$ of equal size,
$\wh{S}_{A,B}(u)$ is the submatrix of $\wh{S}(u)$ on the rows $A$ and the columns $B$, and
$\big|\,\cdot\,\big|_{ik}$ is its quasi-determinant at the entry in the row $i$ and the column
$k$. Replacing the row index $3$ by $4$ turns a quasi-determinant of $S(u)$ with last row $3$
into one of $\wh{S}(u)$. Hence
\aln{
f_{31}(u)\,\wh{d}_1(u) &= \wh{s}_{41}(u) &&= \wh{f}_{41}(u)\,\wh{d}_1(u) ,
\\[0.2em]
f_{32}(u)\,\wh{d}_2(u) &= \big|\,\wh{S}_{\{1,4\},\,\{1,2\}}(u)\,\big|_{42} &&= \wh{f}_{42}(u)\,\wh{d}_2(u) ,
\\[0.2em]
d_3(u) &= \big|\,\wh{S}_{\{1,2,4\},\,\{1,2,3\}}(u)\,\big|_{43} &&= \wh{f}_{43}(u)\,\wh{d}_3(u) .
}
The series $e_{23}(u)$ involves only the rows $1,2$, so
$\wh{d}_2(u)\,e_{23}(u) = \wh{d}_2(u)\,\wh{e}_{23}(u)$. By Remark \ref{R:no-shift} the relations
\eqref{fork:ef+triv} hold with every index shifted by $1$; in particular
$\wh{e}_{23}(u) = \wh{f}_{42}(-u) = b(-u)$ and $\wh{f}_{43}(u) = 0$. This gives \eqref{blocks:2}.
\end{proof}

By Lemma \ref{L:block} the Gaussian decomposition of the block $S^\circ(u)$ reads
\eqa{
s_{11}(u) &= d_1(u) , \qq  & s_{12}(u) &= d_1(u)\,e_{12}(u) , \\[0.2em]
s_{21}(u) &= f(u)\,d_1(u) , \qq & s_{22}(u) &= d_2(u) + f(u)\,d_1(u)\,e_{12}(u) , \\[0.2em]
s_{31}(u) &= g(u)\,d_1(u), \qq & s_{32}(u) &= b(u)\,d_2(u) + g(u)\,d_1(u)\,e_{12}(u) ,
\label{deg:gauss}
}
on the entries that we shall use.

We now collect the relations of the block that we shall use. On the positions $1,2$,
\cite[Lem.~5.1]{Re26} applies at $i=1$ through Lemma \ref{L:AI-block}. This gives \eqref{AI:dd},
\eqref{AI:f=e}, \eqref{AI:df} and their (5.7) at $i=1$. On the positions $2,3$, Corollary
\ref{C:o4o3} and Remark~\ref{R:no-shift} give the first relation of \eqref{fork:bb} for $b(u)$. It
has the form of their (5.7) at $i=2$ and $\ka=2$, with $d_3(u)=0$. Finally
$[d_1(u),b(v)]=0$ by Lemma \ref{L:emb-comm} with $m=1$, because $d_1(u)=\wh{d}_1(u)$ and the
coefficients of $b(v)$ lie in $X^{\tw[1]}$ by Lemma~\ref{L:Shat-red}.

The relations between $f(u)=f_{21}(u)$ and $b(u)=f_{32}(u)$ that we need are the following.

\begin{prop} \label{P:deg-AI}
In the block $S^\circ(u)$ the following relations hold:
\gat{
f(u)\,b(v) = \frac{u-v-1}{u-v}\,b(v)\,f(u)
+ \frac{\big[\,b^{(0)},f(u)\,\big] - \big[\,b(v),f^{(0)}\,\big]}{u-v} ,
\label{deg:ff1-block}
\\[0.3em]
\Sym_{w,v} \big[b(w),[b(v),f(u)]\big] = 0 ,
\label{deg:serre-ii}
\\[0.3em]
\Sym_{w,v} \big[f(w),[f(v),b(u)]\big] = \frac{1}{v+w-1}\,\Sym_{w,v}
\bigg( \frac{(2v-1)\,b(u)}{(u-v)(u+v-1)} - \frac{b(v)}{u-v} + \frac{b(1-v)}{u+v-1} \bigg)
d^{-1}_1(v)\,d_2(v) .
\label{deg:serre-iii}
}
\end{prop}

Relation \eqref{deg:ff1-block} has the form of \cite[(5.8)]{Re26} for the series $f(u)$ and
$b(u)$, and \eqref{deg:serre-ii} and \eqref{deg:serre-iii} have the form of their (5.9)
with the doubled series $b(u)$ and $f(u)$, respectively.
The zero right hand side of \eqref{deg:serre-ii} agrees in form with their (5.9) at
$d_3(v)=0$, but does not follow from the arguments of \cite[\S4]{LWZ25a}, which invert $d_3(v)$.

The proof of Proposition \ref{P:deg-AI} occupies the rest of this section.
It rests on Lemma \ref{L:deg-comm} below, which gives the relations involving $g(u)$ that the proof
needs and the {\it cross} relation \eqref{deg:ff3} between $f(u)$ and $b(v)$. It also uses the
identities \eqref{deg:b0}, \eqref{deg:f0} and \eqref{deg:fb} from the proof of that lemma.


\subsection{The relations involving $g(u)$} \label{sec:deg-g}

We use divided differences. For a series $a(u)$ put
\[
a[x,y] := \frac{a(x)-a(y)}{x-y} \quad\text{ and }\quad a[u,v,w] := \frac{a(u)}{(u-v)(u-w)} + \frac{a(v)}{(v-u)(v-w)} + \frac{a(w)}{(w-u)(w-v)} .
\]
Both are symmetric in their arguments, and $a[u,v,w] = \big(a[u,w]-a[v,w]\big)/(u-v)$. If $a(u)$
has zero constant term, then the coefficient of $u^{-1}$ in $a[u,v]$ is $-a(v)$. Put
$\xi(x) := d^{-1}_1(x)\,d_2(x)$. In this notation \cite[(5.7)]{Re26} at $i=1$ reads
\equ{
[\,f(x),\,f(y)\,] = f[x,y]\,\big(f(y)-f(x)\big) + \frac{\xi(y)-\xi(x)}{x+y-1} ,
\label{deg:AI}
}
and the second relation of \eqref{AI:f=e} at $i=1$ gives $\xi(1-x)=\xi(x)$.

\begin{lemma} \label{L:deg-comm}
In the block $S^\circ(u)$ the following relations hold:
\ali{
f(u)\,g(v) &= g(v)\,f(u) + \frac{\big(g(u)-g(v)\big)\big(f(v)-f(u)\big)}{u-v}
+ \frac{b(v)\,d^{-1}_1(v)\,d_2(v) - b(1-u)\,d^{-1}_1(u)\,d_2(u)}{u+v-1} ,
\label{deg:fg} \\[0.2em]
f(u)\,b(v) &= \frac{b(v)\,\big(f(v)+(u-v-1)\,f(u)\big) + g(u) - g(v)}{u-v} ,
\label{deg:ff3} \\[0.2em]
[\,b(v),\,g(u)\,] &= b(v)\,[\,b(v),\,f(u)\,] = [\,b(v),\,f(u)\,]\,b(v) .
\label{deg:bg}
}
\end{lemma}

\begin{proof}
The argument uses three inner derivations of $X^\tw(\mfo_6,\mfo_2\op\mfo_4)$:
$\mathrm{ad}\,b^{(0)}$, $\mathrm{ad}\,f^{(0)}$ and $\pd := [\,\wh{d}_{1,1},\,\cdot\,]$, where
$\wh{d}_{1,1}$ is the coefficient of $u^{-2}$ in $d_1(u)=\wh{d}_1(u)$. We first compute their
action on the series of the block. The derivation $\mathrm{ad}\,f^{(0)}$ gives the coefficients of $u^{-1}$ in
\eqref{deg:ff3} and \eqref{deg:fg}. The derivation $\pd$ then gives the other coefficients by
induction on the degree in $u^{-1}$. Relation \eqref{deg:bg} follows from \eqref{deg:ff3}.

\smallskip

{\it The derivations.} By \eqref{deg:gauss}, $b^{(0)}$ and $f^{(0)}$ are the coefficients of
$u^{-1}$ in $s_{32}(u)$ and in $s_{21}(u)$. The constant term of $S^\circ(u)$ is
$\diag(1,1,0)$. So the coefficient of $u^{-1}$ in \eqref{RE-iso} gives, for $1\le k,l\le3$,
\eqa{
[\,b^{(0)},\,s_{kl}(v)\,] &= \del_{k2}\,s_{3l}(v) + \del_{l2}\,s_{k3}(v) ,
\\
{}[\,f^{(0)},\,s_{kl}(v)\,] &= \del_{k1}\,s_{2l}(v) - \del_{k2}\,s_{1l}(v) + \del_{l1}\,s_{k2}(v)
- \del_{l2}\,s_{k1}(v) .
\label{deg:s1}
}
We apply \eqref{deg:s1} to the entries \eqref{deg:gauss}. The first relation on the first column
gives \eqref{deg:b0} below. The second relation on the rows $1$ and $2$ gives
\aln{
[\,f^{(0)},d_1(v)\,] &= f(v)\,d_1(v) + d_1(v)\,e_{12}(v) , \qu
[\,f^{(0)},e_{12}(v)\,] = \xi(v) - 1 - e_{12}(v)^2 ,
\\[0.2em]
[\,f^{(0)},d_2(v)\,] &= -\,f(v)\,d_2(v) - d_2(v)\,e_{12}(v) .
}
With these, the row $3$ gives \eqref{deg:f0}:
\gat{
[\,b^{(0)},\,d_1(v)\,] = 0 , \qq [\,b^{(0)},\,f(v)\,] = g(v) , \qq [\,b^{(0)},\,g(v)\,] = 0 ,
\label{deg:b0}
\\[0.2em]
[\,f^{(0)},\,g(v)\,] = b(v)\,\xi(v) - g(v)\,f(v) , \qq
[\,f^{(0)},\,b(v)\,] = b(v)\,f(v) - g(v) .
\label{deg:f0}
}
By \eqref{deg:b0}, $\mathrm{ad}\,b^{(0)}$ commutes with $\pd$, and $g^{(0)} = [b^{(0)},f^{(0)}]$.
Relation \eqref{AI:df} at $i=1$ gives the action of $\pd$ on $f(x)$, and
$\mathrm{ad}\,b^{(0)}$ takes it to the action on $g(x)$. So for $a(x)$ equal to
$f(x)$ or $g(x)$ we have
\gat{
\pd\big(a(x)\big) = (1-2x)\,a(x) + 2a^{(0)} ,
\label{deg:D}
\\[0.2em]
\pd\big(a[u,v]\big) = (1-2u)\,a[u,v] - 2a(v) = (1-u-v)\,a[u,v] - a(u) - a(v) .
\label{deg:Dd}
}
The derivation $\pd$ annihilates $b(x)$ and $\xi(x)$, because $[d_1(u),b(v)]=0$ and the diagonal
series commute.

{\it The induction.} Let $X = \sum_{r,s\ge0} X_{rs}\,u^{-r-1}v^{-s-1}$ be the difference of the
two sides of \eqref{deg:ff3} or of \eqref{deg:fg}. Let $X_0(v)$ and $X^0(u)$ be the coefficients
of $u^{-1}$ and of $v^{-1}$ in $X$. We show that $X_0(v)=0$ and that
\equ{
\pd(X) = (1-2u)\,X + 2X_0(v) , \qq \pd(X) = (2-2u-2v)\,X + 2X_0(v) + 2X^0(u) ,
\label{deg:DX}
}
for \eqref{deg:ff3} and for \eqref{deg:fg}, respectively. The coefficient of $u^{-r-1}v^{-s-1}$
in \eqref{deg:DX} gives $2X_{r+1,s} = X_{rs} - \pd(X_{rs})$ in the first case and
$2X_{r+1,s} = 2X_{rs} - 2X_{r,s+1} - \pd(X_{rs})$ in the second. Since $X_{0s} = 0$ for all $s$,
induction on $r$ gives $X=0$.

\smallskip

{\it Proof of \eqref{deg:ff3}.} The relation reads
\equ{
[\,f(u),\,b(v)\,] = g[u,v] - b(v)\,f[u,v] .
\label{deg:fb}
}
Its coefficient of $u^{-1}$ is the second relation of \eqref{deg:f0}, so $X_0(v)=0$. The first
relation of \eqref{deg:DX} follows from \eqref{deg:D} and \eqref{deg:Dd}, because
$\pd\big(b(v)\big)=0$.

\smallskip

{\it Proof of \eqref{deg:fg}.} Since $\xi(1-u)=\xi(u)$, the relation reads
$[f(u),g(v)] = Y_1+Y_2$ with
\[
Y_1 := (v-u)\,g[u,v]\,f[u,v] , \qq Y_2 := c[v,1-u] , \qq c(x) := b(x)\,\xi(x) .
\]
Since $c(x)$ has zero constant term, $Y_2$ has only negative powers of $u$ and $v$. The coefficient
of $u^{-1}$ in the relation is the first relation of \eqref{deg:f0}, so $X_0(v)=0$. By \eqref{deg:D}
the commutator $[f(u),g(v)]$ satisfies the second relation of \eqref{deg:DX}, with its own
coefficients of $u^{-1}$ and $v^{-1}$ in place of $X_0(v)$ and $X^0(u)$. By \eqref{deg:Dd}, and since
$\pd$ annihilates $b(x)$ and $\xi(x)$,
\[
\pd(Y_1) = (2-2u-2v)\,Y_1 + 2g(u)\,f(u) - 2g(v)\,f(v) , \qq \pd(Y_2) = 0 .
\]
Here $-g(v)\,f(v)$ and $g(u)\,f(u)$ are the coefficients of $u^{-1}$ and of $v^{-1}$ in $Y_1$.
The sum of the coefficients of $u^{-1}$ and of $v^{-1}$ in $Y_2$ is $(u+v-1)\,Y_2$. So $Y_1$ and
$Y_2$ satisfy the second relation of \eqref{deg:DX}, and so does $X$.

\smallskip

{\it Proof of \eqref{deg:bg}.} Apply $\mathrm{ad}\,b^{(0)}$ to \eqref{deg:fb}. We use
\eqref{deg:b0} and $[b^{(0)},b(v)] = -b(v)^2$, which is the coefficient of $u^{-1}$ in the first
relation of \eqref{fork:bb}. This gives
\[
[\,b(v),\,g(u)\,] - [\,b(v)^2,\,f(u)\,] = -\,b(v)\,[\,b(v),\,f(u)\,] ,
\]
that is, $[b(v),g(u)] = [b(v),f(u)]\,b(v)$. So the first equality of \eqref{deg:bg} is equivalent to
$Z(u)=0$, where $Z(x) := [b(v),[b(v),f(x)]]$. The series $Z(x)$ contains only monomials
$x^{-r-1}v^{-s-2}$ with $r,s\ge0$. So each coefficient after putting $v$ for $x$ is a finite sum. Apply $\mathrm{ad}\,b(v)$ to \eqref{deg:fb}, and use
$[b(v),g(x)] = [b(v),f(x)]\,b(v)$ for $x\in\{u,v\}$. This gives $Z(u) = Z[u,v]$, that is,
$(u-v-1)\,Z(u) = -Z(v)$. Let $Z_r$ be the coefficient of $u^{-r-1}$ in $Z(u)$. The coefficients
of $u^0,u^{-1},\ldots$ give $Z_0=-Z(v)$ and $Z_{r+1}=(v+1)\,Z_r$, so $Z_r=-(v+1)^r Z(v)$. Since
$b(v)$ has zero constant term, every $Z_r$ contains only the powers $v^{-2},v^{-3},\ldots$. If
$Z(v)\ne0$ and $v^{-m}$ is its highest power, then $Z_m=-(v+1)^m Z(v)$ has a nonzero term in
$v^0$, which is impossible. Hence $Z(v)=0$, and so $Z(u)=0$.
\end{proof}

Relation \eqref{deg:ff3} links the sub-blocks on the positions $1,2$ and $2,3$. It is the
analogue of the type DI relation \cite[(5.51)]{Re26}.


\subsection{Proof of Proposition \ref{P:deg-AI}} \label{sec:deg-proof}

We obtain \eqref{deg:ff1-block} by eliminating the corner series from \eqref{deg:ff3}. For the
Serre relations we apply $\mathrm{ad}\,b(w)$ or $\mathrm{ad}\,f(w)$ to \eqref{deg:fb}. The first
double commutator is antisymmetric in $v$ and $w$. The second is antisymmetric up to a fraction,
which gives the right hand side of \eqref{deg:serre-iii}.

\begin{proof}[Proof of \eqref{deg:ff1-block}]
By \eqref{deg:b0} and \eqref{deg:f0} we have $g(u) = [b^{(0)},f(u)]$ and
$g(v) = b(v)\,f(v) + [b(v),f^{(0)}]$. Substitute these into \eqref{deg:ff3}. The terms
$b(v)\,f(v)$ cancel, and \eqref{deg:ff1-block} follows.
\end{proof}

\begin{proof}[Proof of \eqref{deg:serre-ii}]
Apply $\mathrm{ad}\,b(w)$ to \eqref{deg:fb}. In the block $S^\circ(u)$ the first relation of \eqref{fork:bb}
reads $[b(w),b(v)] = \big(b(v)-b(w)\big)\,b[v,w]$. Replace $v$ by $w$ in \eqref{deg:fb} and take
the divided difference in $u$. This gives $[b(w),f[u,v]] = b(w)\,f[u,v,w]-g[u,v,w]$. Relation
\eqref{deg:bg} gives $[b(w),g[u,v]] = b(w)\,[b(w),f[u,v]]$. Hence
\equ{
[\,b(w),[\,b(v),f(u)\,]\,]
= \big(b(v)-b(w)\big)\big(b[v,w]\,f[u,v] + b(w)\,f[u,v,w] - g[u,v,w]\big) .
\label{deg:serre-ii-mid}
}
The first factor is antisymmetric in $v$ and $w$. The second factor is symmetric, because
$f[u,v]-f[u,w] = (v-w)\,f[u,v,w]$ and $(v-w)\,b[v,w] = b(v)-b(w)$. So $\Sym_{w,v}$ of
\eqref{deg:serre-ii-mid} vanishes.
\end{proof}

\begin{proof}[Proof of \eqref{deg:serre-iii}]
Interchange $u$ and $v$ in \eqref{deg:fb} and apply $\mathrm{ad}\,f(w)$. For $[f(w),b(u)]$ replace
$(u,v)$ by $(w,u)$ in \eqref{deg:fb}. This gives
\equ{
[\,f(w),[\,f(v),b(u)\,]\,] = [\,f(w),g[u,v]\,] - g[u,w]\,f[u,v]
- b(u)\,\big([\,f(w),f[u,v]\,] - f[u,w]\,f[u,v]\big) .
\label{deg:serre-iii-split}
}
By \eqref{deg:AI} and \eqref{deg:fg} we have
\aln{
[\,f(w),f(x)\,] &= f[x,w]\,\big(f(x)-f(w)\big) + \frac{\xi(x)-\xi(w)}{x+w-1} ,
\\[0.2em]
[\,f(w),g(x)\,] &= g[x,w]\,\big(f(x)-f(w)\big) + \frac{b(x)\,\xi(x)-b(1-w)\,\xi(w)}{x+w-1} .
}
Let $a(x)$ be $f(x)$ or $g(x)$. The divided difference in $x$ of $a[x,w]\,\big(f(x)-f(w)\big)$ is
$a[u,w]\,f[u,v] + a[u,v,w]\,\big(f(v)-f(w)\big)$, and its first summand cancels in
\eqref{deg:serre-iii-split}. Take the divided difference in $x$ at $u,v$ of each fraction. Multiply
the first result by $-b(u)$ on the left. Add the second result. Simplification gives
\equ{
[\,f(w),[\,f(v),b(u)\,]\,] = \big(g[u,v,w]-b(u)\,f[u,v,w]\big)\big(f(v)-f(w)\big)
+ \frac{b[u,v]\,\xi(v)-b[u,1-w]\,\xi(w)}{v+w-1} .
\label{deg:serre-iii-mid}
}
The first term vanishes under $\Sym_{w,v}$, because $g[u,v,w]$ and $f[u,v,w]$ are symmetric in
$v$ and $w$ while $f(v)-f(w)$ is antisymmetric. Under $\Sym_{w,v}$ we may interchange $v$ and $w$
in the part of the second term with $\xi(w)$. Finally
\[
b[u,v]-b[u,1-v] = \frac{(2v-1)\,b(u)}{(u-v)(u+v-1)} - \frac{b(v)}{u-v} + \frac{b(1-v)}{u+v-1} .
\]
This gives \eqref{deg:serre-iii}.
\end{proof}


\section{Drinfeld type relations} \label{sec:Dr}

In this section we prove the Drinfeld type presentation of
$SY^\tw(\mfo_{2n},\mfo_{n-1}\op\mfo_{n+1})$, Theorem \ref{T:iso}. For $n\ge3$ the relations are
verified in \S\ref{sec:rels}. Generation is shown in \S\ref{sec:gen-lem}, and the proof,
including the case $n=2$, is completed in \S\ref{sec:proof-iso}. There we also compare the
presentation with \cite{LZ25} and \cite{HL26}. Section~\ref{sec:ext} treats the extended twisted
Yangian and the fundamental representation, and \S\ref{sec:coideal} the coideal structure.


\subsection{Drinfeld series} \label{sec:Dr-series}

We label the nodes of the Dynkin diagram of $D_n$ as follows. For $n\ge3$ the nodes
$1,\ldots,n-1$ form a chain and the node $n$ is attached to the node $n-2$; for $n=2$ the two
nodes are not connected. The diagram involution $\tau$ fixes the nodes $1,\ldots,n-2$, which we
call the {\it $\iota$-nodes}. It interchanges the two {\it fork nodes} $n-1$ and $n$, which form
the {\it $\tau$-pair}. We write $(c_{ij})_{1\le i,j\le n}$ for the Cartan matrix of $D_n$. %

We attach a pair of Drinfeld series to each node. At an $\iota$-node the spectral parameter is
shifted as in \cite[(5.90)]{Re26}. At the $\tau$-pair there is no shift, by Remark
\ref{R:no-shift}:
\ali{
\label{Dr1}
b_i(u) &:= \wh{f}_{i+1,i}\big(u+\tfrac{\ka-i}{2}\big) , &
h_i(u) &:= \wh{d}^{\,-1}_i\big(u+\tfrac{\ka-i}{2}\big)\,\wh{d}_{i+1}\big(u+\tfrac{\ka-i}{2}\big)
\qu (1\le i\le n-2) ,
\\[0.3em]
\label{Dr2}
b_{n-1}(u) &:= \wh{f}_{n,n-1}(u) , & h_{n-1}(u) &:= \frac{2u}{2u-1}\,\wh{d}^{\,-1}_{n-1}(u)\,\wh{d}_n(u) , & &
\\[0.3em]
\label{Dr3}
b_{n}(u) &:= \wh{f}_{n+1,n-1}(u) , & h_{n}(u) &:= h_{n-1}(-u) . & &
}
The elements $b_{i,r}$ and $h_{i,r}$ are defined by $b_i(u) = \sum_{r\ge0}b_{i,r}u^{-r-1}$ and
$h_i(u) = 1+\sum_{r\ge0}h_{i,r}u^{-r-1}$.

By Remark \ref{R:no-shift}, the relations of Corollary \ref{C:o4o3} hold in
$X^\tw(\mfo_{2n},\mfo_{n-1}\op\mfo_{n+1})$ with every index shifted by $n-2$ and with the same
spectral parameter. This shift turns the series \eqref{o4o3:Dr} into \eqref{Dr2}--\eqref{Dr3}. We
cite the shifted relations by the labels of Corollary \ref{C:o4o3}.

\begin{lemma} \label{L:heven}
$h_i(u) = h_i(-u)$ for $1\le i\le n-2$, so that $h_{i,2r}=0$ for those $i$; whereas
$h_{n-1}(u) = h_n(-u)$, equivalently $h_{n-1,r} = (-1)^{r+1}h_{n,r}$.
\end{lemma}

\begin{proof}
Let $1\le i\le n-2$. By the second relation of \eqref{AI:f=e} and by \eqref{AI:dd} we have
\[
\wh{d}^{\,-1}_i(u)\,\wh{d}_{i+1}(u) = \wh{d}^{\,-1}_i(\ka-i-u)\,\wh{d}_{i+1}(\ka-i-u) .
\]
Replace $u$ by $u+\frac{\ka-i}{2}$. This gives $h_i(u)=h_i(-u)$. The second assertion is the
definition \eqref{Dr3}.
\end{proof}


\subsection{The presentation} \label{sec:presentation}

We now state the main result. Recall the notations $\{f(u)\}^u$ and $\Sym$ from
\S\ref{sec:notation}, and $\{a,b\}=ab+ba$.

\begin{thm} \label{T:iso}
The special twisted Yangian $SY^\tw(\mfo_{2n},\mfo_{n-1}\op\mfo_{n+1})$ is isomorphic to the
associative unital $\C$-algebra generated by the elements $b_{i,r}$ with $1\le i\le n$,
$h_{i,2r+1}$ with $1\le i\le n-2$, and $h_{n-1,r}$, where $r\ge0$, assembled into the series
\eqg{
b_i(u) = \sum_{r\ge0} b_{i,r} \, u^{-r-1}  \;\; (1\le i\le n) , \qq
h_i(u) = 1 + \sum_{r\ge0} h_{i,2r+1} u^{-2r-2} \;\; (1\le i\le n-2) , \\
h_{n-1}(u) = 1 + \sum_{r\ge0} h_{n-1,r} u^{-r-1} ,
\label{iso:series}
}
with $h_n(u) := h_{n-1}(-u)$, and satisfying relations given by the following identities, in
which $1\le i,j\le n$ and rational functions of $u,v,w$ are to be expanded in the domain
$|u|>|v|>|w|$. The Cartan series commute,
\equ{
[h_i(u),\,h_j(v)] = 0 .
\label{iso:hh}
}
If $i$ and $j$ are not both fork nodes, then
\equ{
h_i(u)\,b_j(v) = \bigg( \frac{(u-v+\frac{c_{ij}}{2})(u+v-\frac{c_{ij}}{2})}{(u-v-\frac{c_{ij}}{2})(u+v+\frac{c_{ij}}{2})}\,b_j(v) - \bigg\{ \frac{c_{ij} (u-\tfrac{c_{ij}}2)}{ u\,(u-v-\tfrac{c_{ij}}2)} \, b_j(u-\tfrac{c_{ij}}2) \bigg\}^{\!u}\,  \bigg)\,h_i(u) ,
\label{iso:hb}
}
whereas for $i\in\{n-1,n\}$
\ali{
[h_i(u),\,b_i(v)] &= -\,\frac{\{h_i(u),\,b_i(u)-b_i(v)\}}{u-v} , \qq &
[h_i(u),\,b_{\tau i}(v)] &= \frac{\{h_i(u),\,b_{\tau i}(-u)-b_{\tau i}(v)\}}{u+v} .
\label{iso:hb-fork}
}
The relations among the series $b_i(u)$ are
\ali{
(u+v)\,[b_i(u),\,b_j(v)] &= \del_{\tau i,j}\,\big( h_j(v) - h_i(u) \big) \qu\text{if}\qu c_{ij}=0 ,
\label{iso:bb0}
\\[0.2em]
b_i(u)\,b_j(v) &= \frac{u-v+\tfrac{c_{ij}}{2}}{u-v-\tfrac{c_{ij}}{2}}\,b_j(v)\,b_i(u) + \frac{1}{u-v-\tfrac{c_{ij}}{2}}\,\big([b_{i,0},b_j(v)] - [b_i(u),b_{j,0}]\big) \qu\text{if}\qu c_{ij}=-1 ,
\label{iso:bb1}
\\[0.2em]
[b_i(u),\,b_i(v)] &= \frac{(b_i(v) - b_i(u))^2}{v-u} + \del_{\tau i,i}\,\frac{h_i(v) - h_i(u)}{v+u} ,
\label{iso:bb2}
}
and the Serre type relation reads
\equ{
\Sym_{w,v} \big[b_i(w),[b_i(v),b_j(u)]\big] = \frac{\del_{\tau i,i}}{v+w}\,\Sym_{w,v} x_{ij}(v,u)
\qu\text{if}\qu c_{ij} = -1 ,
\label{iso:S}
}
where
\equ{
x_{ij}(v,u) = \bigg( \frac{2v\,b_j(u)}{(u-v-\frac12)(u+v-\frac12)}
- \frac{b_j(v+\frac12)}{u-v-\frac12} + \frac{b_j(\frac12-v)}{u+v-\frac12} \bigg)\,h_i(v) .
\label{iso:x}
}
\end{thm}

In the split presentation \cite[Thm.~5.14]{Re26} in type D the factors $\del_{\tau i,j}$ and
$\del_{\tau i,i}$ are $\del_{ij}$ and $1$, relation \eqref{iso:hb} holds for all pairs of nodes,
and the Cartan series are even and independent. Here the Cartan term of \eqref{iso:bb2} at a fork
node moves to \eqref{iso:bb0}, which couples the two fork nodes. The right hand side of
\eqref{iso:S} vanishes when the doubled node is a fork node; this reflects $d_3(u)=0$ in Lemma
\ref{L:block}. Relations \eqref{iso:hh}, \eqref{iso:hb-fork}, \eqref{iso:bb0} and
\eqref{iso:bb2} at the fork nodes become \eqref{Y:drinfeld} under
$x^-(u)=b_{n-1}(u)$, $x^+(u)=b_n(-u)$ and $\xi(u)=h_{n-1}(u)$
(Theorem \ref{T:o4o3}).


\subsection{The relations} \label{sec:rels}

In this subsection $n\ge3$. Proposition \ref{P:mixed} treats the pairs formed by the node $n-2$
and a fork node, and Proposition \ref{P:rels} collects the relations of Theorem \ref{T:iso}.

\begin{lemma} \label{L:dd}
We have $[\wh{d}_i(u),\wh{d}_j(v)] = 0$ for all $1\le i,j\le n$.
\end{lemma}

\begin{proof}
For $i,j\le n-1$ this is \eqref{AI:dd}. For $n-1\le i,j\le n$ it follows from the first relation of \eqref{fork:hh}
through Lemma \ref{L:Shat-red}. Otherwise one index is at most $n-2$ and the other is $n$.
By Lemma \ref{L:Shat-red} the coefficients of $\wh{d}_n(v)$ lie in
$X^{\tw[n-2]}$, and the claim is Lemma \ref{L:emb-comm} with $m=n-2$.
\end{proof}

\smallskip

{\it Commutation with the fork series.} By Lemma \ref{L:Shat-red} the coefficients of
$b_{n-1}(v)$ and $b_n(v)$ lie in $X^{\tw[n-2]}$. So Lemma \ref{L:emb-comm} with $m=n-2$ gives
\equ{
[\wh{f}_{j+1,j}(u), b_i(v)] = 0 \qu\text{and}\qu [\wh{d}_k(u), b_i(v)] = 0 \qu\text{for}\qu j<n-2,\; k<n-1,\; i=n-1,n . \label{mixed:triv}
}

\begin{lemma} \label{L:dn-comm}
We have $[\wh{d}_n(u),\, b_j(v)] = 0$ for every $\iota$-node $1\le j\le n-2$.
\end{lemma}

\begin{proof}
By Lemma \ref{L:Shat-red} with $m=n-1$, the coefficients of $\wh{d}_n(u)$
lie in $X^{\tw[n-1]}$. Here $b_j(v) = \wh{f}_{j+1,j}\big(v+\tfrac{\ka-j}{2}\big)$ and $j<n-1$. So
the claim is Lemma \ref{L:emb-comm} with $m=n-1$.
\end{proof}

\begin{prop} \label{P:mixed}
Let $j\in\{n-1,n\}$. Then relations \eqref{iso:hb}, \eqref{iso:bb1} and \eqref{iso:S} hold for the
pairs of nodes $(n-2,j)$ and $(j,n-2)$.
\end{prop}

\begin{proof}
We use \eqref{d-b} and \cite[(5.15)]{Re26} for \eqref{iso:hb}, and Proposition \ref{P:deg-AI}
and the automorphism $\vartheta$ for \eqref{iso:bb1} and \eqref{iso:S}. Since
$\frac{\ka-n+2}{2}=\frac12$, the definition \eqref{Dr1} at $i=n-2$ reads
\equ{
b_{n-2}(u) = \wh{f}_{n-1,n-2}\big(u+\tfrac12\big) , \qq
h_{n-2}(u) = \wh{d}^{\,-1}_{n-2}\big(u+\tfrac12\big)\,\wh{d}_{n-1}\big(u+\tfrac12\big) .
\label{mixed:Dr}
}
Write $(k,l)$ for either pair. Then $c_{kl}=-1$.

\smallskip

{\it Relation \eqref{iso:hb}.} Let $(k,l) = (n-2,j)$ and $u' := u+\tfrac12$. By \eqref{mixed:Dr},
$h_{n-2}(u)\,b_j(v) = \wh{d}^{\,-1}_{n-2}(u')\,\wh{d}_{n-1}(u')\,b_j(v)$. Apply \eqref{d-b} at
$u'$ to the last two factors. By \eqref{mixed:triv} the coefficients of $\wh{d}^{\,-1}_{n-2}(u')$
commute with those of $b_j(v)$. So $\wh{d}^{\,-1}_{n-2}(u')$ moves back to the right, and by
\eqref{mixed:Dr} the factor on the right becomes $h_{n-2}(u)$. The substitution $u\mapsto u'$
turns the three coefficients of \eqref{d-b} into those of \eqref{iso:hb} after expanding
$\{\,\cdot\,\}^u$, and the
arguments $u$ and $1-u$ into $u+\frac12$ and $\frac12-u$.

Let $(k,l) = (j,n-2)$. By \eqref{Dr2} and Lemma \ref{L:dn-comm},
$\wh{d}_n(u)$ commutes past $b_{n-2}(v)$ in $h_{n-1}(u)b_{n-2}(v)$.
By Lemma \ref{L:AI-block}, apply \cite[(5.15)]{Re26} with $i=n-2$, $\ka-i=1$ and
$v$ replaced by $v+\frac12$ to the remaining product
$\wh{d}^{\,-1}_{n-1}(u)\,\wh{f}_{n-1,n-2}(v+\frac12)$.
Using \eqref{mixed:Dr} and \eqref{Dr2} gives \eqref{iso:hb} for $(k,l)=(n-1,n-2)$.
The case $(k,l)=(n,n-2)$ follows on replacing $u$ by $-u$. Indeed $h_n(u) = h_{n-1}(-u)$ by
\eqref{Dr3}, and in \eqref{iso:hb} the rational coefficient of $b_l(v)$ and the bracket
$\{\,\cdot\,\}^u$ are even in $u$.

\smallskip

{\it Relations \eqref{iso:bb1} and \eqref{iso:S}.} First let $j=n$, and take the series
$f(u)$, $b(u)$ of \eqref{deg:notation}. By Lemma \ref{L:reduce-to-3} and \eqref{blocks:1}, the
relations of Proposition \ref{P:deg-AI} hold for every $n\ge3$, with $d_1(u)$ and $d_2(u)$
replaced by $\wh{d}_{n-2}(u)$ and $\wh{d}_{n-1}(u)$. By \eqref{mixed:Dr} and \eqref{Dr3},
$b_{n-2}(u) = f(u+\frac12)$ and $b_n(u) = b(u)$, so $b_{n-2,0} = f^{(0)}$ and $b_{n,0} = b^{(0)}$.

Relation \eqref{iso:bb1} for $(k,l)=(n-2,n)$ is \eqref{deg:ff1-block} with $u$ replaced by
$u+\frac12$. For the reversed pair, multiply \eqref{iso:bb1} by its denominator:
\[
(u-v+\tfrac12)\,b_k(u)\,b_l(v) - (u-v-\tfrac12)\,b_l(v)\,b_k(u) = [b_{k,0},b_l(v)] - [b_k(u),b_{l,0}] .
\]
This identity is unchanged when $(k,u)$ and $(l,v)$ are interchanged. So \eqref{iso:bb1} for
$(n,n-2)$ is equivalent to \eqref{iso:bb1} for $(n-2,n)$.

In \eqref{iso:S} we have $\del_{\tau k,k}=0$ at a fork node and $\del_{\tau k,k}=1$ at the
$\iota$-node $n-2$. So \eqref{iso:S} for $(k,l)=(n,n-2)$ is \eqref{deg:serre-ii} with $u$
replaced by $u+\frac12$. For $(k,l)=(n-2,n)$, replace $v$ by $v+\frac12$ and $w$ by $w+\frac12$ in
\eqref{deg:serre-iii}. By \eqref{mixed:Dr} the last factor becomes $h_{n-2}(v)$, and by
\eqref{iso:x} this gives \eqref{iso:S}.

Now let $j=n-1$. By Lemma \ref{L:sigma} and \eqref{Dr1}--\eqref{Dr3}, the automorphism
$\vartheta$ fixes $b_{n-2}(u)$ and $h_{n-2}(u)$ and maps $b_n(u)$ to $b_{n-1}(u)$. It maps
\eqref{iso:bb1} and \eqref{iso:S} for $j=n$ to the same relations for $j=n-1$.
\end{proof}

\begin{prop} \label{P:rels}
Let $n\ge3$. The Drinfeld series \eqref{Dr1}--\eqref{Dr3} satisfy the relations
\eqref{iso:hh}--\eqref{iso:x} of Theorem \ref{T:iso} in
$X^\tw(\mfo_{2n},\mfo_{n-1}\op\mfo_{n+1})$, and hence in $SY^\tw(\mfo_{2n},\mfo_{n-1}\op\mfo_{n+1})$.
\end{prop}

\begin{proof}
We sort the pairs of nodes $(i,j)$ by the number of fork nodes among them. Relation
\eqref{iso:hh} is Lemma~\ref{L:dd} for all $i,j$.

Let $1\le i,j\le n-2$. By Lemma \ref{L:AI-block} the type AI relations \eqref{AI:dd} and
\cite[(5.5)--(5.10)]{Re26} hold for these $i,j$, with the same constant $\ka$. In the Drinfeld
normalisation \eqref{Dr1} they are \cite[(5.106)--(5.110)]{Re26}, as rewritten in the proof of
their Thm.~5.14. Here $\del_{\tau i,i}=1$ and $\del_{\tau i,j}=\del_{ij}$, so these are the
relations \eqref{iso:hb}, \eqref{iso:bb0}--\eqref{iso:bb2} and \eqref{iso:S} for those $i,j$.

Let one of $i,j$ be a fork node and the other be $n-2$. Then $c_{ij}=-1$, and
\eqref{iso:hb}, \eqref{iso:bb1} and \eqref{iso:S} are Proposition \ref{P:mixed}.

Let one of $i,j$ be a fork node and the other be at most $n-3$. Then $c_{ij}=0$ and
$\del_{\tau i,j}=0$. So \eqref{iso:hb} reads $[h_i(u),b_j(v)]=0$ and \eqref{iso:bb0} reads
$[b_i(u),b_j(v)]=0$. The first holds by \eqref{mixed:triv} if $j$ is the fork node, and by
\eqref{AI:dd} and Lemma \ref{L:dn-comm} if $i$ is the fork node. The second holds by
\eqref{mixed:triv}.

Let $i,j\in\{n-1,n\}$. Relation \eqref{iso:hb-fork} at $i=n-1$ is \eqref{fork:hb}; at $i=n$ it
follows from \eqref{fork:hb} on replacing $u$ by $-u$, since $h_n(u)=h_{n-1}(-u)$. Relation
\eqref{iso:bb2} has $\del_{\tau i,i}=0$ and is the first relation of \eqref{fork:bb}. Relation
\eqref{iso:bb0} for $(i,j)=(n-1,n)$ is the second relation of \eqref{fork:bb}, and for
$(i,j)=(n,n-1)$ it follows on interchanging $u$ and $v$. \qedhere
\end{proof}

\begin{remark} \label{R:AIII}
The base case $n=3$ of Lemma \ref{L:reduce-to-3} admits a second description: under
$\mfo_6\cong\mfsl_4$ the pair $(\mfo_6,\mfo_2\op\mfo_4)$ becomes
$(\mfsl_4,\mathfrak{s}(\mfgl_2\op\mfgl_2))$, the quasi-split pair of type AIII$_3$, with the
$\iota$-node $1$ and the $\tau$-pair $\{2,3\}$ becoming the middle node and the two end nodes of
$A_3$. This case is treated by a Gaussian decomposition in \cite[\S7.3]{LZ25}. Several
relations appear there as identities of generating series, among them their (7.19), the
analogue of \eqref{deg:ff3}, and the Serre relations with a doubled end node in their
Prop.~7.14. The Serre relation with the doubled middle node and the remaining Drinfeld relations
are recorded modewise in their Props.~7.14 and 7.15. In quasi-split type A, \cite[App.~B]{HL26} expresses the modes $h_{i,r}$ and $b_{i,r}$ with $r=0,1$
through the Drinfeld generators of the Yangian. Their appendix does not derive closed identities
among the full series $h_i(u)$ and $b_i(u)$. Proposition \ref{P:rels} at $n=3$ gives these identities
in type AIII$_3$ in the form of Theorem \ref{T:iso}.
\end{remark}


\subsection{Generation} \label{sec:gen-lem}

We show that the coefficients of the Drinfeld series generate the special twisted Yangian. For
$n\ge3$ their symbols are, up to an automorphism of $\mfo_{2n}[x]^\rho$, the generators of
\cite[Prop.~2.5]{LZ25}.

\smallskip

{\it The involution $\rho$.} Let $\pi$ be the involution of $\{1,\ldots,N\}$ with
$\pi(i)=\ol{\si(i)}$. It interchanges $i$ and $\ol{\imath}$ for $i\ne n,n+1$, fixes $n$ and
$n+1$, and commutes with $a\mapsto\bar{a}$. We have $X^t = -JXJ$ for $X\in\mfo_{2n}$, since
$F^t_{ij} = -J\,F_{ij}\,J$. With \eqref{G-props} this gives
$\rho(X) = -G\,X^t\,G^{-1} = G\,J\,X\,J\,G = (GJ)\,X\,(GJ)^{-1}$, and $GJ = \sum_i E_{i\pi(i)}$ is
the permutation matrix of $\pi$ by \eqref{G}. Since $\pi$ commutes with $a\mapsto\bar{a}$, we obtain
\equ{
\rho(F_{ij}) = F_{\pi(i)\pi(j)} .
\label{gen:rho}
}

\smallskip

{\it The generators of \cite{LZ25}.} The Chevalley generators $e_i,f_i$ of $\mfo_{2n}$ are
$f_i = F_{i+1,i}$ and $e_i = F_{i,i+1}$ for $i\le n-1$, and $f_n = F_{n+1,n-1}$, $e_n = F_{n-1,n+1}$.
Put $H_i := [e_i,f_i]$. By \eqref{gen:rho}, $\rho(e_i)=-f_{\tau(i)}$ and $\rho(f_i)=-e_{\tau(i)}$. So
$\rho$ is the involution $\omega_0$ of \cite[\S2.2]{LZ25}, and $\mfo_{2n}[x]^\rho$ is their twisted
current algebra $\mf{g}[z]^{\check\omega}$, with $z=x$. Let $n\ge3$. By their Prop.~2.5 and its
proof, the elements
\equ{
x_{i,r} := \big((-1)^rf_i-e_{\tau(i)}\big)\,x^r , \qq
\xi_{i,r} := \big(H_{\tau(i)}-(-1)^rH_i\big)\,x^r \qq (1\le i\le n,\; r\ge0)
\label{gen:LZ}
}
of $\mfo_{2n}[x]^\rho$ generate $U(\mfo_{2n}[x]^\rho)$. Here $\xi_{i,2r}=0$ for $i\le n-2$, and
$\xi_{n,r}=(-1)^{r+1}\xi_{n-1,r}$.

The symbols computed below are not the elements \eqref{gen:LZ}. They differ from them by the
exchange of the fork nodes and by the sign $(-1)^r$ in degree $r$. We now define the automorphism
that records this.

Let $P_\si$ be the permutation matrix of $\si$. Since $\si$ commutes with $a\mapsto\bar{a}$,
conjugation by $P_\si$ maps $F_{ij}$ to $F_{\si(i)\si(j)}$. So it maps $e_i$, $f_i$ and $H_i$ to
the same elements at the node $\tau(i)$. It commutes with $\rho$ by \eqref{gen:rho}, since
$\si$ commutes with $\pi$. Hence the automorphism
\[
\varsigma : X x^r\mapsto(-1)^r\,P_\si XP_\si\,x^r
\]
of $\mfo_{2n}[x]$ commutes with the involution $Xx^r\mapsto(-1)^r\rho(X)\,x^r$, and it restricts to
an automorphism of the graded Lie algebra $\mfo_{2n}[x]^\rho$. It satisfies
$\varsigma(x_{i,r})=(-1)^rx_{\tau(i),r}$ and $\varsigma(\xi_{i,r})=(-1)^r\xi_{\tau(i),r}$.

\smallskip

{\it Symbols of the Drinfeld series.} In the filtration of \S\ref{sec:notation} the
coefficients $b_{i,r}$ and $h_{i,r}$ lie in $\mr{F}_r$. As in the proof of Lemma
\ref{L:o4o3-symbols}, every nonlinear Gaussian correction has degree at most $r-1$;
the same bound holds for the terms with at least two factors in \eqref{Dr1} and \eqref{Dr2}.
Since $\wh{s}_{ij}(u) = s_{\si(i)j}(u)$,
the symbol of $\wh{f}^{\,(r)}_{ji}$ is $\bar{s}^{(r+1)}_{\si(j)i}$ and that of $\wh{d}_{i,r}$ is
$\bar{s}^{(r+1)}_{\si(i)i}$. The shifts of the spectral parameter in \eqref{Dr1} and the prefactor
$\frac{2u}{2u-1}$ in \eqref{Dr2} add only terms of lower degree, except the scalar $\tfrac12$ that
the prefactor adds to $h_{n-1,0}$. The scalars
$\ga^{(r)}_{ij} = \tfrac12\del_{r0}(g_{ij}-\del_{i\ol{\jmath}})$ of \eqref{SY-graded-iso} occur
only in degree $0$. They vanish for $b_{i,0}$, since the entries that occur are $(i+1,i)$ with $i\le n-2$,
$(n+1,n-1)$ and $(n,n-1)$, and for these $g_{ab}=\del_{a\ol{b}}=0$. For $h_{n-1,0}$ they give
$\ga^{(0)}_{n+1,n}-\ga^{(0)}_{n-1,n-1}=-\tfrac12$, which cancels that scalar. By
\eqref{SY-graded-iso} we obtain
\eqg{
\bar{b}_{i,r} \mapsto F^{(\rho,r)}_{i+1,i} \;\; (1\le i\le n-2) , \qq
\bar{b}_{n-1,r} \mapsto F^{(\rho,r)}_{n+1,n-1} , \qq
\bar{b}_{n,r} \mapsto F^{(\rho,r)}_{n,n-1} ,
\\[0.3em]
\bar{h}_{i,2r+1} \mapsto F^{(\rho,2r+1)}_{i+1,i+1} - F^{(\rho,2r+1)}_{ii} \;\; (1\le i\le n-2) ,
\qq
\bar{h}_{n-1,r} \mapsto F^{(\rho,r)}_{n+1,n} - F^{(\rho,r)}_{n-1,n-1} .
\label{gen:symbols}
}
Since $g_{aj}=\del_{a\si(j)}$, \eqref{SY-graded-iso} gives
$F^{(\rho,r)}_{ij} = \big(F_{i\si(j)}-(-1)^rF_{j\si(i)}\big)\,x^r$. For example,
\gan{
F^{(\rho,r)}_{n+1,n-1} = \big(f_n-(-1)^re_{n-1}\big)\,x^r = (-1)^rx_{n,r} ,
\\
F^{(\rho,r)}_{n+1,n}-F^{(\rho,r)}_{n-1,n-1} = \big((-1)^rH_{n-1}-H_n\big)\,x^r = (-1)^r\xi_{n,r} .
}
The other rows of \eqref{gen:symbols} are computed in the same way. By Lemma \ref{L:heven} the
coefficients $h_{i,2r}$ with $i\le n-2$ and $h_{n,r}$, absent from \eqref{gen:symbols}, satisfy the
relations stated for $\xi_{i,2r}$ and $\xi_{n,r}$ after \eqref{gen:LZ}. Hence, for all $i$ and $r$,
\equ{
\bar{b}_{i,r} \mapsto \varsigma(x_{i,r}) , \qq \bar{h}_{i,r} \mapsto \varsigma(\xi_{i,r}) .
\label{gen:varsigma}
}

\begin{lemma} \label{L:gen}
The elements $b_{i,r}$ with $1\le i\le n$ and $h_{i,2r+1}$ with $1\le i\le n-2$, together with
$h_{n-1,r}$, $r\ge0$, generate the algebra $SY^\tw(\mfo_{2n},\mfo_{n-1}\op\mfo_{n+1})$.
\end{lemma}

\begin{proof}
For $n=2$ this is shown in the proof of Theorem \ref{T:o4o3}. Let $n\ge3$. By
\eqref{gen:varsigma}, the symbols of the listed elements are the images under $\varsigma$ of the
elements \eqref{gen:LZ} other than $\xi_{i,2r}$ with $i\le n-2$ and $\xi_{n,r}$. These omitted
elements are $0$ and multiples of $\xi_{n-1,r}$. Since $\varsigma$ is an automorphism, the symbols
generate $U(\mfo_{2n}[x]^\rho)\cong\gr SY^\tw(\mfo_{2n},\mfo_{n-1}\op\mfo_{n+1})$. An induction on
the filtration degree, as in the proof of Theorem \ref{T:o4o3}, gives the lemma.
\end{proof}


\subsection{Proof of the presentation} \label{sec:proof-iso}

We first deduce mode relations from the current relations of Theorem \ref{T:iso}. We then prove
injectivity from these relations and the classical presentation of \cite[Prop.~2.7]{LZ25}. We use the
finite type Serre relations of their Thm.~3.9 in the form \cite[Def.~2.10]{HL26}.

\begin{lemma} \label{L:LZ}
Let $n\ge3$, and let $b_i(u)$ and $h_i(u)$ be series in a unital associative $\C$-algebra,
of the form \eqref{iso:series}, with $h_n(u) = h_{n-1}(-u)$, and satisfying
\eqref{iso:hh}--\eqref{iso:x}.
For $1\le i\le n$ and $r\ge0$, denote their coefficients by $b_{i,r}$ and $h_{i,r}$.
Set $h_{i,-1}=1$ and $h_{i,r}=0$ for $r<-1$. Then these elements satisfy the following relations, in which $r,s\ge0$, except
that $r\ge-1$ in \eqref{LZ:hb}:
\ali{
[h_{i,r},h_{j,s}] &= 0 , \qq h_{\tau i,r} = (-1)^{r+1} h_{i,r} ,
\label{LZ:hh}
\\[0.2em]
[h_{i,r+1},b_{j,s}] - [h_{i,r-1},b_{j,s+2}] &=
\tfrac{1}{2}(c_{ij}-c_{\tau i,j})\,\{h_{i,r},b_{j,s}\}
+ \tfrac{1}{2}(c_{ij}+c_{\tau i,j})\,\{h_{i,r-1},b_{j,s+1}\}
\nn\\
& \qu + \tfrac{1}{4}c_{ij}c_{\tau i,j}\,[h_{i,r-1},b_{j,s}] ,
\label{LZ:hb}
\\[0.2em]
[b_{i,r+1},b_{j,s}] - [b_{i,r},b_{j,s+1}] &= \tfrac{1}{2}c_{ij}\,\{b_{i,r},b_{j,s}\}
- 2\,\del_{\tau i,j}\,(-1)^{r}\,h_{j,r+s+1} ,
\label{LZ:bb}
}
together with the Serre type relations
\ali{
[b_{i,r},b_{j,s}] &= \del_{\tau i,j}\,(-1)^r\,h_{j,r+s} &&\text{if } c_{ij}=0 ,
\label{LZ:S0}
\\
\big[b_{i,0},[b_{i,0},b_{j,0}]\big] &= -\,\del_{\tau i,i}\,b_{j,0} \qq &&\text{if } c_{ij}=-1 .
\label{LZ:S1}
}
These are the relations (3.1)--(3.6) and (3.10) of \cite{LZ25} at $\hbar=1$, and their finite
type Serre relations (3.18) and (3.19). The cases $r=-1$ and $r=0$ of \eqref{LZ:hb} are their
(3.2) and (3.3). The case $r=s=0$ of \eqref{LZ:S0} is their finite type Serre
relation (3.17) together with the case $j=\tau i$, which is omitted in
their Thm.~3.9. Relations \eqref{LZ:hh}--\eqref{LZ:bb}, \eqref{LZ:S1} and the case $r=s=0$
of \eqref{LZ:S0} contain the defining relations of $\mc{Y}^\imath$ in type {\rm DII}$_n$ in the
form \cite[Def.~2.10]{HL26}.
\end{lemma}

\begin{proof}
Relations \eqref{LZ:hh} are \eqref{iso:hh} and \eqref{iso:series}. Relation \eqref{LZ:hb} follows
by clearing the denominators in \eqref{iso:hb} and \eqref{iso:hb-fork} and comparing the
coefficients of $u^{-r}v^{-s-1}$; cf.~\cite[(3.22)]{LZ25}. After clearing the denominators, the
terms with $b_j(\pm u-\frac{c_{ij}}{2})$ or $b_j(\pm u)$ have no negative power of $v$. The
conventions on $h_{i,r}$ record the constant term of $h_i(u)$, and $c_{\tau i,j}=c_{ij}$ unless $i$
and $j$ are both fork nodes.

Divide \eqref{iso:bb0} by $u+v$. Then the Cartan terms of \eqref{iso:bb0} and \eqref{iso:bb2} are
$\del_{\tau i,j}(h_j(v)-h_i(u))/(u+v)$, where $j=i$ in \eqref{iso:bb2}. If $j=\tau i$, then
$h_i(u)=h_j(-u)$ by \eqref{iso:series}, and
\[
\frac{h_j(v)-h_j(-u)}{u+v} = \sum_{r,s\ge0}(-1)^r h_{j,r+s}\,u^{-r-1}v^{-s-1} .
\]
This identity turns \eqref{iso:bb0} into \eqref{LZ:S0}. If $c_{ij}=0$, then \eqref{LZ:bb} is the
difference of the cases $(r+1,s)$ and $(r,s+1)$ of \eqref{LZ:S0}. Otherwise \eqref{LZ:bb} follows by
clearing the denominators in \eqref{iso:bb1} and \eqref{iso:bb2} and comparing coefficients.

The coefficient of $v^{-1}w^{-1}u^{-1}$ in \eqref{iso:S} is twice \eqref{LZ:S1}. On the right
hand side only $x_{ij}(w,u)/(v+w)$ contributes, and the coefficient of $w^{-1}u^{-1}$ in
$x_{ij}(w,u)$ is $-2b_{j,0}$.
\end{proof}

\begin{proof}[Proof of Theorem \ref{T:iso}]
Let $\mr{Y}$ denote the associative unital $\C$-algebra with the generators and relations listed
in the theorem. For $n=2$ there are no $\iota$-nodes. In the currents $b^-(u)=b_1(u)$ and
$b^+(u)=b_2(-u)$ the relations of the theorem become the Drinfeld presentation
\eqref{Y:drinfeld} with $x^\pm=b^\pm$ and $\xi=h_1$, so the statement is Theorem \ref{T:o4o3}. Let
$n\ge3$. Then $\mfo_{2n}$ is simple, and the results of \cite{LZ25} quoted below apply. By Lemma~\ref{L:heven} and Proposition \ref{P:rels} the assignment
\equ{
\Phi : \mr{Y} \longrightarrow SY^\tw(\mfo_{2n},\mfo_{n-1}\op\mfo_{n+1}) , \qu
b_{i,r} \mapsto b_{i,r} , \qu h_{i,r} \mapsto h_{i,r} ,
\label{Phi}
}
respects the defining relations of $\mr{Y}$. So it extends to a homomorphism of algebras, which is
surjective by Lemma \ref{L:gen}. We prove that $\Phi$ is injective by passing to the associated
graded algebras, using the classical presentation \cite[Prop.~2.7]{LZ25}.

\smallskip

{\it Step 1: a surjection onto $\gr\mr{Y}$.} Equip $\mr{Y}$ with the filtration
$\deg b_{i,r}=\deg h_{i,r}=r$, and write $\bar{b}_{i,r}$ and $\bar{h}_{i,r}$ for the symbols of
these elements in $\gr\mr{Y}$. By \cite[Prop.~2.7]{LZ25} the algebra $U(\mfo_{2n}[x]^\rho)$ is
generated by the elements \eqref{gen:LZ}, subject to the following relations: \eqref{LZ:hh},
\eqref{LZ:S1}, the case $j=\tau i$ of \eqref{LZ:S0} and its case $r=s=0$, in which $h_{i,r}$ and
$b_{i,r}$ are replaced by $\xi_{i,r}$ and $x_{i,r}$; and
\[
[\xi_{i,r},x_{j,s}] = \big(c_{ij}-(-1)^rc_{\tau i,j}\big)\,x_{j,r+s} , \qq
[x_{i,r+1},x_{j,s}]-[x_{i,r},x_{j,s+1}] = -2\del_{\tau i,j}(-1)^r\xi_{j,r+s+1} .
\]
Lemma \ref{L:LZ} applies to $\mr{Y}$, so its relations hold in
$\mr{Y}$. We show that the symbols $\bar{h}_{i,r}$ and $\bar{b}_{i,r}$ satisfy the relations just
listed. Relations \eqref{LZ:hh}, \eqref{LZ:S0} and \eqref{LZ:S1} are homogeneous, and the second
relation above is the leading term of \eqref{LZ:bb}. The leading terms of the cases
$r=-1$, $r=0$ and $r\ge1$ of \eqref{LZ:hb} are
\gan{
[\bar{h}_{i,0},\bar{b}_{j,s}] = (c_{ij}-c_{\tau i,j})\,\bar{b}_{j,s} , \qq
[\bar{h}_{i,1},\bar{b}_{j,s}] = (c_{ij}+c_{\tau i,j})\,\bar{b}_{j,s+1} ,
\\
{}[\bar{h}_{i,r+1},\bar{b}_{j,s}] = [\bar{h}_{i,r-1},\bar{b}_{j,s+2}] ,
}
and they give the first relation above by induction on $r$. So the assignment
$\xi_{i,r}\mapsto\bar{h}_{i,r}$, $x_{i,r}\mapsto\bar{b}_{i,r}$ defines a homomorphism
$\varpi : U(\mfo_{2n}[x]^\rho)\to\gr\mr{Y}$. It is surjective, since $\gr\mr{Y}$ is generated by
the symbols of the generators of $\mr{Y}$.

\smallskip

{\it Step 2: injectivity.} By \eqref{gen:symbols} the homomorphism $\Phi$ is filtered. By
\eqref{gen:varsigma}, the homomorphisms $\gr\Phi\circ\varpi$ and $U(\varsigma)$ agree on the
generators \eqref{gen:LZ} of $U(\mfo_{2n}[x]^\rho)$. So $\gr\Phi\circ\varpi=U(\varsigma)$ is an
isomorphism. Since $\varpi$ is surjective, $\varpi$ and $\gr\Phi$ are isomorphisms. Both
filtrations are exhaustive and bounded below, so $\Phi$ is injective. Hence $\Phi$ is an
isomorphism.
\end{proof}

\begin{crl} \label{C:LWZ}
Let $n\ge3$, and let $\mc{Y}^\imath$ be the twisted Yangian of quasi-split type {\rm DII}$_n$ in
the form \cite[Def.~2.10]{HL26}, with generators $b^\imath_{i,r}$ and $h^\imath_{i,r}$, where
$1\le i\le n$ and $r\ge0$. The assignment
\equ{
\Psi : \mc{Y}^\imath \longrightarrow \mr{Y} , \qu
b^\imath_{i,r} \mapsto b_{i,r} , \qu h^\imath_{i,r} \mapsto h_{i,r}
\label{psi-hom}
}
extends to an isomorphism of algebras. Hence $\Xi := \Phi\circ\Psi$ is an isomorphism, and
\equ{
\mc{Y}^\imath \;\cong\; \mr{Y} \;\cong\; SY^\tw(\mfo_{2n},\mfo_{n-1}\op\mfo_{n+1}) .
\label{iso-chain}
}
In particular the $R$-matrix presentation of \cite{GR16} and the Drinfeld presentation of
\cite{LZ25} agree in quasi-split type DII$_n$.
\end{crl}

\begin{proof}
By Lemma \ref{L:LZ}, the defining relations of $\mc{Y}^\imath$ hold in $\mr{Y}$, so
\eqref{psi-hom} extends to a homomorphism. It is surjective, since its image contains the
generators of $\mr{Y}$. For injectivity we repeat Steps 1 and 2 of the proof of Theorem
\ref{T:iso} for $\mc{Y}^\imath$ and $\Psi$.

\smallskip

{\it The relations used in Step 1.} In this paragraph $b_{i,r}$ and $h_{i,r}$ stand for
$b^\imath_{i,r}$ and $h^\imath_{i,r}$. The defining relations of $\mc{Y}^\imath$ contain the relations
used in Step 1, except the second relation of \eqref{LZ:hh} for $r\ge1$ and \eqref{LZ:S0} beyond the
case $r=s=0$. The former holds by \cite[(3.10)]{LZ25}: it follows from its case $r=0$ and
\eqref{LZ:bb}. We now prove \eqref{LZ:S0} for all $r,s\ge0$. Let $c_{ij}=0$ and put $D_{r,s} := [b_{i,r},b_{j,s}]-\del_{\tau i,j}(-1)^rh_{j,r+s}$.
Relation \eqref{LZ:bb} gives $D_{r+1,s} = D_{r,s+1}$, so $D_{r,s} = D_{0,r+s}$; and $D_{0,0}=0$ is
a defining relation of $\mc{Y}^\imath$. Put
$\tl{h}_{i,1} := h_{i,1}-\tfrac12h^2_{i,0}$. The cases $r=-1$ and $r=0$ of \eqref{LZ:hb} give
$[\tl{h}_{i,1},b_{k,s}] = (c_{ik}+c_{\tau i,k})\,b_{k,s+1}$ for all $k$ and $s$. With
\eqref{LZ:hh} and $D_{1,t}=D_{0,t+1}$ this gives
\[
[\tl{h}_{i,1},D_{0,t}] = (c_{ii}+c_{\tau i,i}+c_{ij}+c_{\tau i,j})\,D_{0,t+1} .
\]
If $j=\tau i$, then $c_{ii}+c_{\tau i,i} = c_{ij}+c_{\tau i,j} = 2$, and the terms with
$h_{j,t+1}$ cancel. If $j\ne\tau i$, then $c_{ij}+c_{\tau i,j}=0$. In both cases
$c_{ii}+c_{\tau i,i}\ge2$, so the scalar factor is positive. By induction on $t$
every $D_{0,t}$ vanishes.

\smallskip

{\it Injectivity of $\Psi$.} Equip $\mc{Y}^\imath$ with the filtration
$\deg b^\imath_{i,r} = \deg h^\imath_{i,r} = r$. The relations of $\mc{Y}^\imath$ now contain
those used in Step 1, so Step 1 gives a surjective homomorphism
$\varpi^\imath : U(\mfo_{2n}[x]^\rho)\to\gr\mc{Y}^\imath$ with $\xi_{i,r}\mapsto\bar{h}^\imath_{i,r}$
and $x_{i,r}\mapsto\bar{b}^\imath_{i,r}$. The homomorphism $\Psi$ is filtered, and
$\gr\Psi\circ\varpi^\imath = \varpi$. The map $\varpi$ is an isomorphism by Step 2 of the proof of
Theorem \ref{T:iso}. Hence $\varpi^\imath$ is injective, and so $\varpi^\imath$ and $\gr\Psi$ are
isomorphisms. Therefore $\Psi$ is injective, and it is an isomorphism. In particular
$\varpi^\imath$ is the isomorphism $U(\mfo_{2n}[x]^\rho)\cong\gr\mc{Y}^\imath$ of
\cite[Cor.~3.8]{LZ25} in type DII$_n$; we do not use that result.

\smallskip

{\it The last assertion.} The twisted Yangian of \cite{GR16} is isomorphic to
$SY^\tw(\mfo_{2n},\mfo_{n-1}\op\mfo_{n+1})$ by \cite[Thm.~4.1, Cor.~5.2]{GR16}.
The relation $[b_{i,0},b_{\tau i,0}]=h_{\tau i,0}$ is included in \cite[Def.~2.10]{HL26} and in
\cite[Def.~3.1]{LZ25}, but it must be added in their Thm.~3.9.
The preceding $D_{r,s}$ argument supplies its all-mode consequence, so the proof of their
Thm.~3.9 identifies $\mc{Y}^\imath$ with the algebra of their Def.~3.1 at $\hbar=1$. Here their
Serre relation (3.8) is read with $h_{i,-1}=1$, as in \cite[Thm.~5.1, Rem.~5.2]{LWZ25a} and
\cite[Lem.~5.15]{Re26}. With the convention $h_{i,s}=0$ for $s<0$ stated at the end of
\cite[Def.~3.1]{LZ25}, the right hand side of their (3.8) at $k_1=k_2=r=0$ would vanish,
whereas their (3.19) gives $[b_{i,0},[b_{i,0},b_{j,0}]]=-b_{j,0}$. This identification is
the only place where their (3.8) enters.
\end{proof}

\begin{crl} \label{C:J}
Let $n\ge3$. Then $SY^\tw(\mfo_{2n},\mfo_{n-1}\op\mfo_{n+1})$ is isomorphic to the twisted
Yangian of the pair \eqref{our-pair} in the $J$ presentation of \cite{Mac02}, in the form
\cite[Def.~2.7]{HL26}.
\end{crl}

\begin{proof}
For $n\ge4$ the pair \eqref{our-pair} is of type DI$_n$ in \cite[Tab.~1]{HL26}, and for $n=3$
it is of type AIII$_3$ by Remark \ref{R:AIII}. The underlying simple Lie algebras have types
$\mathsf{D}_n$ and $\mathsf{A}_3$, so neither has type $\mathsf{A}_{2m}$. Hence, by
their Thm.~3.4(3), the algebra $\mc{Y}^\imath$ in type DII$_n$ is isomorphic to the twisted
Yangian in the $J$ presentation. Corollary \ref{C:LWZ} gives the claim.
\end{proof}

\begin{remark} \label{R:Lu}
Corollary \ref{C:J} can alternatively be obtained from the classification of twisted Yangians as
coideal subalgebras of the Yangian by Lu \cite{Lu26b}.
The embedding \eqref{S->TGT} realises $SY^\tw(\mfo_{2n},\mfo_{n-1}\op\mfo_{n+1})$ as a left
coideal subalgebra, and by \cite[Thm.~3.16]{GRW19} this transfers to the $J$ presentation. The
automorphism $X\mapsto X$, $J(X)\mapsto-J(X)$ of the Yangian reverses the coproduct and preserves
the classical limit, so it gives an isomorphic right coideal subalgebra. By
\cite[Thm.~A, Rem.~2.9]{Lu26b} this algebra is isomorphic to the twisted Yangian of
\cite[Def.~2.7]{HL26}. With their Thm.~3.4(3) this also gives an isomorphism
$\mc{Y}^\imath\cong SY^\tw(\mfo_{2n},\mfo_{n-1}\op\mfo_{n+1})$, but not in terms of Drinfeld
generators.
\end{remark}


\subsection{The extended twisted Yangian} \label{sec:ext}

Write $\mc{P} := \C[\ms{t}_0,\ms{t}_1,\ldots]$ for
the algebra of polynomials in countably many variables, and keep the notation $b_{i,r}$,
$h_{i,r}$ for the coefficients of the Drinfeld series \eqref{Dr1}--\eqref{Dr3}, which are
defined in $X^\tw(\mfo_{2n},\mfo_{n-1}\op\mfo_{n+1})$ already.

\begin{crl} \label{C:iso-ext}
Let $n\ge3$. The assignment
\equ{
\Theta \;:\; \mc{Y}^\imath\ot\mc{P} \longrightarrow
X^\tw(\mfo_{2n},\mfo_{n-1}\op\mfo_{n+1}) , \qu
b^\imath_{i,r}\ot1 \mapsto b_{i,r} , \qu h^\imath_{i,r}\ot1 \mapsto h_{i,r} ,
\label{Omega}
}
$1\ot\ms{t}_{2r}\mapsto\ms{c}_{2r}$ and $1\ot\ms{t}_{2r+1}\mapsto\ms{z}_{2r+1}$, defines an
isomorphism of algebras. Equivalently, by Corollary \ref{C:LWZ},
$X^\tw(\mfo_{2n},\mfo_{n-1}\op\mfo_{n+1})$ is presented by the generators and relations of
Theorem \ref{T:iso} together with countably many central generators
$\ms{c}_{2r},\ms{z}_{2r+1}$, $r\ge0$. This refines \eqref{X=SYcz}, whose two polynomial factors
are realised as the images of $\C[\ms{t}_{2r}]$ and $\C[\ms{t}_{2r+1}]$; it is the quasi-split
analogue of \cite[Rem.~5.18]{Re26}.
\end{crl}

\begin{proof}
By normalisation invariance and Lemma \ref{L:X-vs-SY}, the coefficients $b_{i,r},h_{i,r}$
lie in the subalgebra $\mc{S}$ of Section~\ref{sec:norm}; the quotient isomorphism
$\mc{S}\cong SY^\tw$ sends them to the corresponding coefficients there, so Corollary
\ref{C:LWZ} identifies $\mc{Y}^\imath$ with $\mc{S}$ on the stated generators.
Proposition \ref{P:centre} identifies $\mc{P}$ with the centre $\mc{Z}$ by the displayed
assignment of the $\ms{t}_r$.
Under these identifications, $\Theta$ is the multiplication isomorphism
$\mc{S}\ot\mc{Z}\to X^\tw$ of Section~\ref{sec:norm}.
\end{proof}

\begin{remark}[Fundamental representation] \label{R:rep}
The fundamental representation of \cite[App.~A]{Re26} applies to \eqref{G} unchanged, and
$\ms{c}(u)$ acts by the identity. Formulas (A.2), (A.4), (A.6) therein apply
with $\vka:=q-p=2$ in place of $N-2n$. In (A.2) replace the prefactor $[u/(u-\frac14)]^{\vka}$ by
$\frac{2u}{2u-1}$; this gives $\wh{d}_i(u)$ for $1\le i\le n-1$. Formulas (A.4) and (A.6) give
$h_i(u)$ and $b_i(u)$ for $1\le i\le n-2$. The two innermost
diagonal series carry no prefactor and are given by \eqref{rep-dn} instead. For the $\tau$-pair series,
\allowdisplaybreaks
\gat{
h_n(u)\mapsto h_{n-1}(-u), \qu b_n(u)\mapsto G\,b_{n-1}(u)\,G, \qu \wh{d}_{n+1}(u)\mapsto G\,\wh{d}_n(u)\,G,
\\[0.2em]
\wh{d}_n(u) \mapsto I - \sum_{k=1}^{n-1} \bigg( \frac{E_{kk}}{(u-\frac\ka2)^2} + \frac{E_{\ol{k}\ol{k}}}{(u+\frac\ka2)^2}
- \frac{\vka\,u\,E_{k\ol{k}}}{(u-\frac\ka2)^2(u+\frac\ka2)^2} \bigg)
- \frac{(2u-1)\,E_{nn} - (2u+1)\,E_{n+1,n+1}}{(u-\frac\ka2)(u+\frac\ka2)} ,
\label{rep-dn}
\\[0.5em]
b_{n-1}(u) \mapsto \frac{E_{n+1,n-1}-E_{n+2,n}}{u+\frac\ka2}
- \frac{(u+\frac\ka2-1)\,E_{n-1,n} - (u+\frac\ka2+1)\,E_{n+1,n+2}}{(u-\frac\ka2)(u+\frac\ka2)} ,
\label{rep-bn}
\\[0.5em]
h_{n-1}(u) \mapsto I + \frac{(\ka-1)E_{n-1,n-1} - (2u-1)E_{nn} + (2u+1)E_{n+1,n+1}
- (\ka+1)E_{n+2,n+2} - \vka\,E_{n-1,n+2}}{(u-\frac\ka2)(u+\frac\ka2)} .
\label{rep-hn}
}
\end{remark}


\subsection{Coideal structure} \label{sec:coideal}

Recall from \cite[\S5.6]{Re26} that $SY^\tw(\mfo_{2n},\mfo_{n-1}\op\mfo_{n+1})$ is a left coideal
subalgebra of $Y(\mfo_{2n})$ with the coaction
\equ{
\Delta\big(s_{ij}(u)\big) = \sum_{1\le a,b\le N} t_{ia}\big(u-\tfrac\ka2\big)\,
t_{jb}\big(\tfrac\ka2-u\big) \ot s_{ab}(u) \;\in\;
Y(\mfo_{2n})\ot SY^\tw(\mfo_{2n},\mfo_{n-1}\op\mfo_{n+1})[[u^{-1}]] ,
\label{coid-Ds}
}
where the $t_{ij}(u)$ are the (normalised) generating series of $Y(\mfo_{2n})$.

Let $Y(\mfo_{2n})_{Q_+}$ be the span of the elements of strictly positive root degree, where
$\deg t_{ij}(u) = \eps_i-\eps_j$, $\eps_1,\ldots,\eps_N$ are the weights of the vector
representation and $\eps_{\ol{\imath}} = -\eps_i$. Let $e_i,f_i$ be the Chevalley generators of
$\mfo_{2n}$ labelled as in Section~\ref{sec:gen-lem}. We write
$t_{ij}(u) = \del_{ij}+\sum_{r\ge1}t^{(r)}_{ij}u^{-r}$, so that $t^{(1)}_{ij} = F_{ij}$, and
$\ms{s}_{ab} := s^{(1)}_{ab}-s^{(1)}_{ba}$. For an unordered pair $\{l,m\}$ we write
$t^{(1)}_{[l,m]}$ for the entry whose row index is the smaller one.

\begin{prop} \label{P:coideal}
We have
\ali{
\Delta(b_{i,0}) &= \big(f_{\tau(i)}-e_i\big)\ot1 + 1\ot b_{i,0} \qq (1\le i\le n) ,
\label{coid-b0}
\\[0.2em]
\Delta(h_{n-1,0}) &= -\,2\,F_{nn}\ot1 + 1\ot h_{n-1,0} ,
\label{coid-h0}
\\[0.2em]
\Delta(h_{i,1}) &= \ms{h}_{i,1}\ot1 + 1\ot h_{i,1} + \Omega_i \qq (1\le i\le n-2) ,
\label{coid-h1}
\\[0.2em]
\Delta(h_{n-1,1}) &= \ms{h}_{n-1,1}\ot1 + 1\ot h_{n-1,1} - 2\,F_{nn}\ot h_{n-1,0}
+ \Omega_{n-1} ,
\label{coid-h1-fork}
}
where $\ms{h}_{i,1}\in Y(\mfo_{2n})$ and
$\Omega_i\in Y(\mfo_{2n})_{Q_+}\ot SY^\tw(\mfo_{2n},\mfo_{n-1}\op\mfo_{n+1})$ is given by the
following formulas, in which $l$ and $a$ run over $1,\ldots,N$ unless restricted:
\gat{
\Omega_i = \sum_{l\ne i+1} t^{(1)}_{[l,\,i+1]}\ot\ms{s}_{l,i+1}
\;-\; \sum_{l\ne i} t^{(1)}_{[l,\,i]}\ot\ms{s}_{l,i} \qq (1\le i\le n-2) ,
\label{coid-omega}
\\[0.2em]
\Omega_{n-1} = \sum_{a<n}\big( t^{(1)}_{a,n}\ot s^{(1)}_{an} - t^{(1)}_{a,n+1}\ot s^{(1)}_{n+1,a}
+ t^{(1)}_{n+1,\bar{a}}\ot s^{(1)}_{\bar{a}n} - t^{(1)}_{n,\bar{a}}\ot s^{(1)}_{n+1,\bar{a}}\big)
- \sum_{a\ne n-1} t^{(1)}_{[a,\,n-1]}\ot\ms{s}_{a,n-1} .
\label{coid-omega-fork}
}
\end{prop}

\begin{proof}
We expand \eqref{coid-Ds} to the second order in $u^{-1}$. At the nodes $i\le n-2$ the
calculation of \cite[Prop.~5.20]{Re26} applies. The $\tau$-pair needs the entries of $G$ off the
diagonal.

By \eqref{G}, $g_{ij} = \del_{i\si(j)}$. So the coefficient of $u^{-1}$ in \eqref{coid-Ds} is
\equ{
\Delta\big(s^{(1)}_{ij}\big) = \big(t^{(1)}_{i\si(j)} - t^{(1)}_{j\si(i)}\big)\ot1
+ 1\ot s^{(1)}_{ij} \qq (1\le i,j\le N) .
\label{coid-Ds1}
}
The corresponding formula for diagonal $G$ is \cite[(5.142)]{Re26}. The coefficient of $u^{-2}$ is
\cite[(5.143)]{Re26}; there $G$ enters only the term $\ms{q}_{ij}\ot1$, with
$\ms{q}_{ij}\in Y(\mfo_{2n})$. For an element $x$ of the special twisted Yangian we write
$\text{cross}\,\Delta(x)$ for $\Delta(x)-1\ot x$ modulo $Y(\mfo_{2n})\ot1$. If
$\Delta(x) = x'\ot1+1\ot x$ and $\Delta(y) = y'\ot1+1\ot y$, then
$\text{cross}\,\Delta(xy) = x'\ot y+y'\ot x$.

{\it The first modes.} Since $\wh{S}(\infty)=I$, the quasi-determinantal formulas give
$\wh{f}^{\,(0)}_{ji} = \wh{s}^{\,(1)}_{ji}$ and $\wh{d}_{k,0} = \wh{s}^{\,(1)}_{kk}$. The
argument shift in \eqref{Dr1} does not affect the first mode. By
$\wh{s}_{ij}(u) = s_{\si(i)j}(u)$ we obtain
\equ{
b_{i,0} = s^{(1)}_{i+1,i} \;\; (1\le i\le n-2) , \qq
b_{n-1,0} = s^{(1)}_{n+1,n-1} , \qq b_{n,0} = s^{(1)}_{n,n-1} ,
\label{coid-b0-elt}
}
and $\wh{d}_{k,0} = s^{(1)}_{kk}$ for $k<n$, $\wh{d}_{n,0} = s^{(1)}_{n+1,n}$. By
\eqref{coid-Ds1}, the first legs of the coproducts \eqref{coid-b0-elt} are $f_i-e_i$,
$f_n-e_{n-1}$ and $f_{n-1}-e_n$ in the labelling of Section~\ref{sec:gen-lem}, that is,
$f_{\tau(i)}-e_i$. This is \eqref{coid-b0}. Since
$t^{(1)}_{n+1,n+1}-t^{(1)}_{nn} = F_{n+1,n+1}-F_{nn} = -2F_{nn}$, relation \eqref{coid-Ds1} also
gives
\equ{
\Delta\big(\wh{d}_{k,0}\big) = 1\ot\wh{d}_{k,0} \qu (1\le k\le n-1) , \qq
\Delta\big(\wh{d}_{n,0}\big) = -2F_{nn}\ot1+1\ot\wh{d}_{n,0} .
\label{coid-d0}
}

{\it The nodes $i\le n-2$.} Put $D_k := \text{cross}\,\Delta\big(\wh{d}_{k,1}\big)$. Let
$1\le k\le n-1$. The leading $k\times k$ blocks of $\wh{S}(u)$ and $S(u)$ coincide, and
$g_{ll}=1$ for $l\le k$. So \eqref{coid-Ds1} and the calculation in the proof of
\cite[Prop.~5.20]{Re26} give
\equ{
D_k = \sum_{a\ne k} t^{(1)}_{[a,k]}\ot\ms{s}_{ak} \qq (1\le k\le n-1) .
\label{coid-Dk}
}
Let $i\le n-2$. By Lemma \ref{L:heven}, $h_{i,0}=0$, so the argument shift in \eqref{Dr1} does
not affect $h_{i,1}$, and $h_{i,1}$ is given by \cite[(5.144)]{Re26} in the modes of
$\wh{d}_i(u)$ and $\wh{d}_{i+1}(u)$. Its products of the modes $\wh{d}_{k,0}$ have no cross part
by \eqref{coid-d0}. Hence $\text{cross}\,\Delta(h_{i,1}) = D_{i+1}-D_i$, which gives
\eqref{coid-h1} and \eqref{coid-omega}.

{\it The $\tau$-pair.} By \eqref{Dr2},
$\big(1-\tfrac12u^{-1}\big)\,h_{n-1}(u) = \wh{d}^{\,-1}_{n-1}(u)\,\wh{d}_n(u)$. Hence
\equ{
h_{n-1,0} = \tfrac12 + \wh{d}_{n,0} - \wh{d}_{n-1,0} , \qq
h_{n-1,1} = \tfrac12\,h_{n-1,0} + \wh{d}_{n,1} - \wh{d}_{n-1,1} + \wh{d}^{\,2}_{n-1,0}
- \wh{d}_{n-1,0}\,\wh{d}_{n,0} .
\label{coid-h-elt}
}
The first identity and \eqref{coid-d0} give \eqref{coid-h0}. Next,
$\wh{d}_{n,1} = s^{(2)}_{n+1,n} - \sum_{a<n} s^{(1)}_{n+1,a}\,s^{(1)}_{an}$. Since $\ol{n}=n+1$,
we have $F_{n+1,n}=F_{n,n+1}=0$, that is, $t^{(1)}_{n+1,n} = t^{(1)}_{n,n+1} = 0$. So
\cite[(5.143)]{Re26} and \eqref{coid-Ds1} give
\equ{
D_n = \sum_{a<n}\big( t^{(1)}_{an}\ot s^{(1)}_{an} - t^{(1)}_{a,n+1}\ot s^{(1)}_{n+1,a}
+ t^{(1)}_{n+1,\bar{a}}\ot s^{(1)}_{\bar{a}n} - t^{(1)}_{n\bar{a}}\ot s^{(1)}_{n+1,\bar{a}}\big)
- 2F_{nn}\ot\wh{d}_{n,0} .
\label{coid-Dn}
}
In \eqref{coid-h-elt} the term $\frac12h_{n-1,0}$ has no cross part by \eqref{coid-h0}. By
\eqref{coid-d0}, the term $-\wh{d}_{n-1,0}\,\wh{d}_{n,0}$ contributes $2F_{nn}\ot\wh{d}_{n-1,0}$,
and no other product of modes has a cross part. Together with the last term of \eqref{coid-Dn}
this gives $-2F_{nn}\ot h_{n-1,0}$ modulo $Y(\mfo_{2n})\ot1$, by the first identity of
\eqref{coid-h-elt}. Subtracting \eqref{coid-Dk} at $k=n-1$ gives \eqref{coid-h1-fork} and
\eqref{coid-omega-fork}.

Finally, every ambient factor in \eqref{coid-omega} and \eqref{coid-omega-fork} has row index
smaller than its column index, hence positive root degree. So
$\Omega_i\in Y(\mfo_{2n})_{Q_+}\ot SY^\tw(\mfo_{2n},\mfo_{n-1}\op\mfo_{n+1})$.
\end{proof}

The algebra $SY^\tw(\mfo_{2n},\mfo_{n-1}\op\mfo_{n+1})$ quantises the left Lie coideal structure of
\cite[Thm.~3.3]{GR16}. Proposition \ref{P:coideal} differs from its split counterpart
\cite[Prop.~5.20]{Re26} in three ways. First, the first leg of \eqref{coid-b0} is
$f_{\tau(i)}-e_i$, that is, $f_{\tau(i)}$ modulo $Y(\mfo_{2n})_{Q_+}$. At $i=n-1$ this replaces the
first leg $f_{n-1}$ of their (5.139) by $f_n$. At $i\le n-2$ it is unchanged. Second, the $\tau$-pair has the Cartan mode $h_{n-1,0}$, which is
primitive by \eqref{coid-h0}; in the split case $h_{i,0}=0$, since every $h_i(u)$ is even \cite[Lem.~5.11]{Re26}. Third, $\Delta(h_{n-1,1})$ contains
$-2F_{nn}\ot h_{n-1,0}$, whose first leg is a Cartan element.

The modified mode $\wt{h}_{n-1,1} := h_{n-1,1} - \tfrac12h^2_{n-1,0}$ of \cite[\S2.5]{HL26} removes
this term. Indeed, by \eqref{coid-h0}, $\tfrac12\Delta\big(h^2_{n-1,0}\big)
= 2F^2_{nn}\ot1 + 1\ot\tfrac12h^2_{n-1,0} - 2F_{nn}\ot h_{n-1,0}$, and \eqref{coid-h1-fork} gives
\equ{
\Delta\big(\wt{h}_{n-1,1}\big) = \big(\ms{h}_{n-1,1}-2F^2_{nn}\big)\ot1
+ 1\ot\wt{h}_{n-1,1} + \Omega_{n-1} .
\label{coid-h1-tilde}
}
Thus $\wt{h}_{n-1,1}$, like $h_{i,1}$ with $i\le n-2$, is primitive modulo
$Y(\mfo_{2n})_{Q_+}\ot SY^\tw$. This is the analogue of \cite[(3.32)]{HL26}; there the
tensor legs are interchanged because the twisted Yangian is a right coideal subalgebra.


\medskip
{\noindent \bf Acknowledgements.}
The author thanks Kang Lu for suggesting that the approach of \cite{Re26} be extended to quasi-split
type DI, and for correspondence during the preparation of \cite{HL26}. In the course of that
correspondence, an early version of the relations of Section~\ref{sec:low} was shared with the
authors of \cite{HL26}.

This research has been carried out in the framework of the ``Universities' Excellence Initiative''
programme by the Ministry of Education, Science and Sports of the Republic of Lithuania under the
agreement with the Research Council of Lithuania (project No. S-A-UEI-23-6).
AI tools were used in preparing this manuscript. The author takes full responsibility for its content.



\end{document}